\documentclass[10pt,twoside,notitlepage]{amsart} 
\usepackage[french,english]{babel}
\usepackage{amsmath}
\usepackage{amsthm}
\usepackage{amsfonts}
\usepackage{amssymb,amscd,epsf,verbatim,mathtools,framed}
\usepackage{mathrsfs}
\usepackage{graphicx}
\usepackage{latexsym}
\usepackage{lscape}
\usepackage[colorlinks=true]{hyperref}
\hypersetup{colorlinks, citecolor=blue, filecolor=black, linkcolor=purple, urlcolor=violet}
\usepackage{epstopdf}
\usepackage{tikz}
\usetikzlibrary{calc}
\usetikzlibrary{matrix,arrows,decorations.pathmorphing}
\usepackage{tikz-cd}
\usepackage{color}
\usepackage{multirow}  
\usepackage{scalefnt}
\usepackage{fancyhdr}
\newcommand{\Z}{\mathbb{Z}}
\newcommand{\F}{\mathbb{F}}
\newcommand{\N}{\mathbb{N}}
\newcommand{\R}{\mathbb{R}}
\newcommand{\Q}{\mathbb{Q}}

\renewcommand{\P}{\mathbb{P}}
\newcommand{\C}{\mathbb{C}}
\newcommand{\T}{\mathbb{T}}
\newcommand{\G}{\mathbb{G}}

\newcommand{\mA}{\mathcal{A}}

\newcommand{\mC}{\mathcal{C}}
\newcommand{\mD}{\mathcal{D}}
\newcommand{\mE}{\mathcal{E}}
\newcommand{\mF}{\mathcal{F}}

\newcommand{\mH}{\mathcal{H}}

\newcommand{\mR}{\mathcal{R}}

\newcommand{\mT}{\mathcal{T}}

\newcommand{\mW}{\mathcal{W}}
\newcommand{\mX}{\mathcal{X}}
\newcommand{\mY}{\mathcal{Y}}

\newcommand{\fm}{\mathfrak{m}}

\newcommand{\tu}{\textup}

\newcommand{\cl}{\overline}

\newcommand{\ra}{\rightarrow}

\newcommand{\sq}{\widetilde}

\DeclareMathOperator{\im}{im}
\DeclareMathOperator{\Hom}{Hom}

\DeclareMathOperator{\Div}{Div}

\DeclareMathOperator{\GL}{GL}
\DeclareMathOperator{\SL}{SL}
\DeclareMathOperator{\PSL}{PSL}

\newcommand{\surjlra}{\relbar\joinrel\twoheadrightarrow}
\DeclareMathOperator{\lra}{\: \longrightarrow \:}

\newcommand{\gr}[1]{\langle {#1} \rangle}

\DeclareMathOperator{\isom}{\;\xrightarrow{\: {}_{\sim} \:} \;}

\DeclareMathOperator{\Symb}{Symb}

\newcommand{\bi}{\begin{itemize}}
\newcommand{\ei}{\end{itemize}}
\newcommand{\bt}{\begin{theorem}}
\newcommand{\et}{\end{theorem}}
\newcommand{\bbt}{\begin{theorem*}}
\newcommand{\eet}{\end{theorem*}}
\newcommand{\bp}{\begin{proposition}}
\newcommand{\ep}{\end{proposition}}
\newcommand{\bl}{\begin{lemma}}
\newcommand{\el}{\end{lemma}}
\newcommand{\bbl}{\begin{lemma*}}
\newcommand{\eel}{\end{lemma*}}
\newcommand{\bc}{\begin{corollary}}
\newcommand{\ec}{\end{corollary}}
\newcommand{\beg}{\begin{example}}
\newcommand{\eeg}{\end{example}}
\newcommand{\br}{\begin{remark}}
\newcommand{\er}{\end{remark}}
\newcommand{\bbr}{\begin{remark*}}
\newcommand{\eer}{\end{remark*}}
\newcommand{\bd}{\begin{definition}}
\newcommand{\ed}{\end{definition}}
\newcommand{\be}{\begin{enumerate}}
\newcommand{\ee}{\end{enumerate}}
\newcommand{\bex}{\begin{exercise}}
\newcommand{\eex}{\end{exercise}}
\newcommand{\bproof}{\begin{proof}}
\newcommand{\eproof}{\end{proof}}

\theoremstyle{theorem}
\newtheorem{theorem}{Theorem}[section]
 
\newtheorem{assumption}[theorem]{Assumption}

\theoremstyle{definition}
\newtheorem{example}[theorem]{Example}

\newtheorem{question}[theorem]{Question}

\newtheorem{definition}[theorem]{Definition}
\newtheorem{proposition}[theorem]{Proposition}

\newtheorem{lemma}[theorem]{Lemma}
\newtheorem{corollary}[theorem]{Corollary}

\newtheorem{notation}[theorem]{Notation}
\newtheorem{conjecture}[theorem]{Conjecture}
\newtheorem{remark}[theorem]{Remark}

\usepackage{etoolbox}
\patchcmd{\section}{\scshape}{\bfseries}{}{}
\makeatletter
\renewcommand{\@secnumfont}{\bfseries}
\makeatother

\numberwithin{equation}{section}
 
\title[The Halo Conjecture for $\GL_2$]{\Large $\textup{The Halo Conjecture for }\GL_2$}
\author{Hansheng Diao}
\email{hdiao@mail.tsinghua.edu.cn}
\address{Yau Mathematical Sciences Center, Tsinghua University}
\author{Zijian Yao}
\email{zijianyao@uchicago.edu}
\address{Department of Mathematics, University of Chicago}

\begin{document}

\maketitle

\begin{abstract}
We prove the Halo conjecture on the geometry of  the eigencurve over the boundary of the weight space, predicted by  Coleman--Mazur and Buzzard--Kilford. 
\end{abstract}

\section{Introduction}
In this paper we study the geometry of the eigencurve over $\Q$ over the boundary of the weight space. For the setup, let $p$ be an odd prime.\footnote{We exclude the case $p = 2$ only to simplify notations (also note that the case $p = 2$ is already treated by \cite{Buzzard_Kilford}).} Let $v(\cdot)$ and $|\cdot|$ denote the $p$-adic valuation and the $p$-adic norm normalized so that $v(p) = 1$ and $|p|=\frac{1}{p}$. Let $\mW$ denote the \emph{weight space} for $\GL_2$, namely the rigid analytic space associated to the Iwasawa algebra $\Lambda = \Z_p [\![\Z_p^\times]\!]$. For convenience, we fix an isomorphism $(\Z_p^\times)_{\tu{tor}} \cong \F_p^\times$ and let $\mW_{\omega}$ denote the component of $\mW$ indexed by a character $\omega: \F_p^\times \ra \Z_p^\times$. For a closed point $\chi:\Z_p^\times \ra \C_p^\times$ in $\mW$, we define $T_{\chi} := \chi(\exp (p)) -1.$ This gives us a parameter on $\mW$. For any $r \in (0,1)$, let $\mW^{> r} = \{|T| > r\}$ denote the \emph{Halo} of the weight space (of radius $r$), which is a union of annuli, and similarly define $\mW^{> r}_{\omega}$ for each component. Now let us fix a tame level 
and let $\mX$ denote the corresponding eigencurve constructed by Coleman--Mazur and Buzzard \cite{Coleman_Mazur, Buzzard_eigen}. Recall that $\mX$ parameterizes (normalized) finite slope overconvergent eigenforms with the fixed  tame level  and it admits a \emph{weight map} $\tu{wt}: \mX \ra \mW$ as well as a \emph{slope map} $a_p: \mX \ra \G_m^{\tu{rig}}.$
\[
\begin{tikzcd}
\mX \arrow[r, "a_p"] \arrow[d, "\tu{wt}"] 
& \G_m^{\tu{rig}} \\
\mW
\end{tikzcd}
\]
At classical points in $\mX$, these maps send a normalized classical eigenform $f$ of weight $k+2$ and $p$-Nebentypus character $\epsilon: (\Z/p^m \Z)^\times \ra \C_p^\times$ to its \textit{weight character} 
\[\tu{wt}(f) = (k, \epsilon): x \mapsto x^k\cdot \epsilon (x)\] and its $U_p$-eigenvalue, respectively. For $r \in (0, 1)$, let $\mX^{> r} := \tu{wt}^{-1} (\mW^{> r})$, and similarly define $\mX^{> r}_{\omega}$ for each $\omega$. While the geometry of $\mX$ over the center of the weight space seems  mysterious, over the boundary the structure of $\mX$ is rather well-behaved. In particular, we have the following \emph{Halo conjecture}, suggested by the work Buzzard--Kilford \cite{Buzzard_Kilford} and Coleman--Mazur \cite{Coleman_Mazur} (see also the thesis of Emerton \cite{Emerton} for a related computation). 
\begin{conjecture}[The Halo conjecture, {\cite[Conjecture 1.2]{LWX}}] \label{conj:Halo}
There exists some $r \in (0, 1)$ sufficiently close to $1$ such that the following holds. 
\be
\item $\mX^{> r}$ is a disjoint union of countably infinitely many components $\mX_{1}^{> r}, \mX_{2}^{> r}, ...$ such that the weight map $\tu{wt}: \mX_n^{> r} \ra \mW^{> r}$ is finite flat for each $n$.
\item There are rational numbers $\alpha_1 \le \alpha_2 \le  ... \in \Q_{\ge 0}$ tending to infinity, such that for each component $\mX_{n}^{> r}$ and each point $z \in \mX_{n}^{> r}$,  we have \[
v (a_p (z)) = \alpha_n \cdot v(T_{\tu{wt}(z)}). 
\]
\item The sequence $\alpha_1, \alpha_2, ...$ (counted with multiplicities) is a disjoint union of finitely many arithmetic progressions. 
\ee
\end{conjecture}

In its crude form, the main result of this paper is that 
\bt 
\label{mainthm:halo}
The Halo Conjecture (Conjecture \ref{conj:Halo}) holds.  
\et 

We refer the reader to Theorem \ref{mainthm:disjoint_union} and Theorem \ref{mainthm:arithmic_progression} for a more refined version of Theorem \ref{mainthm:halo}. 
Let us make some remarks on the content of our theorem as well as some literature concerning this conjecture.
\br 
\be 
\item As mentioned above, Buzzard--Kilford proved this conjecture for $p = 2$, using explicit coordinates coming from the geometry of the modular curve with small level. Their result played an important role in the recent study of symmetric power functoriality by Newton--Thorne \cite{Newton_Thorne}. 
\item The first significant breakthrough towards the Halo conjecture was made by Liu--Wan--Xiao \cite{LWX} using an ingenious argument for a definite quarternion algebra $D$ over $\Q$. In particular, using the ($p$-adic) Jacquet--Langlands correspondence for $\GL_2$, their work gives part of Conjecture \ref{conj:Halo}, namely the connected components of $\mX$ away from the ones whose tame parts are all principal series. However, the geometry of the eigencurve with tame level $1$ is inaccessible using their approach. 
\item In a slightly more recent work, Ren--Zhao \cite{Ren_Zhao} obtain partial parts on Part (1) of Conjecture \ref{conj:Halo}. Their method is to consider a real quardratic extension $F$ of $\Q$ in which $p$ splits, and then prove a suitable version of Theorem \ref{mainthm:disjoint_union} for a definite quarternion algebra over $F$. At the end of their introduction, it is remarked that more ideas seem to be needed to prove the full Halo conjecture.   
\item The method considered in this paper has suitable generalizations to symplectic groups (for example $\tu{GSp}_4$), which we will investigate in subsequent works.
\ee 
\er 

Next we briefly explain the idea of the proof. Our main strategy is to apply the method of Liu--Wan--Xiao \cite{LWX} for modular symbols. Let $\Delta^0 := \Div^0 (\P^1(\Q))$ be the abelian group of degree zero divisors on $\P^1(\Q)$, which admits a left action of $\GL_2(\Q)$ via linear fractional transformations. Recall that, for a congruence subgroup $\Gamma \subset \PSL_2(\Z)$ and an abelian group $V$ endowed with a right action of $\Gamma$, the \textit{$V$-valued modular symbols} (of level $\Gamma$) are elements of the abelian group 
\begin{equation} \label{eq:maindef_Manin_relation}
\Symb_{\Gamma}(V) = \Hom_{\Gamma} (\Delta^0, V)
\end{equation}
consisting of $\Gamma$-invariant maps from $\Delta^0$ to $V$. One of the difficulties in carrying out the arguments from \cite{LWX} to the modular symbols is the following:  while the ``$V$-valued'' automorphic forms considered in \textit{loc.cit.} are essentially finitely many copies of $V$ (since they consider a definite quarternion algebra which is compact at the infinite place), the 
modular symbols are generally \emph{not} of this form. 

In \cite{Pollack_Stevens}, Pollack--Stevens proved a structural result on the $\Gamma$-action of $\Delta^0$. It is convenient to adopt the following notation: recall that an element $\gamma = \big(\begin{smallmatrix} a & b \\ c & d \end{smallmatrix} \big) \in \textup{PGL}_2 (\Z)$ represents an (oriented) geodesic path $e_\gamma$ in the upper half plane $\mH$,  starting from the point $\frac{b}{d}$ to the point $\frac{a}{c}$ in $\P^1(\Q)$. Its boundary divisor \[[e_{\gamma}] = \partial e_\gamma := \{\frac{a}{c}\} - \{\frac{b}{d}\} \in \Delta^0\] is called a \textit{modular path}. Also recall that such a path $e_\gamma$ is called a \textit{unimodular path} if $\gamma \in \PSL_2(\Z)$. Now suppose $\Gamma$ contains no $2$- or $3$-torsion elements. Then Pollack--Stevens showed that, there exist unimodular paths $e_0, e_1, ..., e_{t-1}, e_\infty$, such that $\Delta^0$ is generated 
by $[e_0], ...,[e_t], [e_\infty]$ as a $\Z[\Gamma]$-module, with a single relation called ``Manin's relation'', generated by
\begin{equation} \label{eq:Manin_relation}
\sum (\gamma_{i}^{-1} - 1) [e_i] = (\gamma_\infty^{-1} - 1) [e_\infty]
\end{equation} 
where $e_\infty = \{\infty\} - \{0\}$, $\gamma_\infty = \big(\begin{smallmatrix} 1 & 1 \\ 0 & 1 \end{smallmatrix} \big)$, and $\gamma_i$ are  elements in $\Gamma$ determined by the choices of $e_i$'s. The presence of Manin's relation (\ref{eq:Manin_relation}) is what prevents us from directly applying the method of \cite{LWX}. However, the result of Pollack--Stevens suggests that the desired freeness of modular symbols might hold if we \emph{remove} some edges from the generating set  $e_i$'s above. This  leads us to consider the space of partial modular symbols. In other words, we will replace $\Delta^0$ by the set of degree $0$ divisors supported on a certain $\Gamma$-invariant subset $C \subset \P^1 (\Q)$. It turns out the eigencurve $\mX_{C}$ associated with such partial modular symbols contains (two copies of) the cuspidal part $\mX^{\tu{cusp}}$ of $\mX$ as a closed subspace. In this paper, we will prove variants of Theorem \ref{mainthm:disjoint_union} and \ref{mainthm:arithmic_progression} for the eigencurve $\mX_{C}$ of parital modular symbols and deduce the halo conjecture from these variants. 

We end the introduction with two standard consequences of the Halo conjecture. Our  first corollary answers a question/conjecture in \cite{Coleman_Mazur}. 

\bc \label{cor:Chenevier}
If $\mC$ is an irreducible component of the eigencurve $\mX$ that is finite over the weight space, then $\mC$ is contained in the ordinary locus. 
\ec

\bc 
Let $k \ge 2$ be an integer. Then every irreducible component of $\mX$ contains a classical point of weight $k$. 
\ec

For the proof of these corollaries, we refer the reader to \cite[Proposition 3.24 and Corollary 1.9]{LWX}, where they prove the corresponding results for the eigencurve of definite quarternion algebras. 
The arguments there  apply  \emph{verbatim} in our setup.
 
\subsection*{Acknowledgement} 
First of all, as it would be clear to the reader, this article owes much debt to the work of Liu--Wan--Xiao \cite{LWX}, which provides some of the key $p$-adic analysis that we need. Additionally, we are very grateful to Rob Pollack, for his encouragement and for pointing out the reference \cite{BD} on partial modular symbols. We also wish to thank Matt Emerton for helpful and clarifying conversations, and for explaining the proof of Buzzard--Kilford \cite{Buzzard_Kilford}.  
 
\section{Partial modular symbols} 
In this section, we briefly recall from \cite{Pollack_Stevens, BD} the notions of  modular symbols and partial modular symbols. The key result we need later is Proposition \ref{prop:partial_free_for_l_square}. 

\subsection{Basic definitions}
Let $\Gamma $ be a finite index subgroup of $\PSL_2 (\Z)$. Let $\Delta:= \Div (\P^1(\Q))$ (resp.  $\Delta^0 := \Div^0 (\P^1(\Q))$) be the abelian group of divisors (resp. degree zero divisors) on $\P^1(\Q)$, which admits a left action of $\GL_2(\Q)$, and thus $\Gamma$, via linear fractional transformations. Let $C \subset \P^1(\Q)$ be a non-empty subset   invariant under the  action of  $\Gamma$. Let
\[
\Delta_C = \{ \sum_{c \in C} n_c \cdot c \: |  \: n_c \in \Z, \,\,n_c = 0 \tu{ for almost all } c \}
\]
be the subgroup of $\Delta$ of divisors supported on $C$, and let $\Delta_C^0$ be its subgroup of degree zero divisors supported on $C$. Note that $\Gamma$ naturally acts on $\Delta_C$ and $\Delta_C^0$. 

\bd[\cite{Pollack_Stevens, BD}] Let $V$ be an abelian group endowed with a right $\Gamma$-action. A $V$-\emph{valued partial modular symbol} (supported on $C$) is a $\Gamma$-invariant homomorphism $\varphi \in \Hom(\Delta^0_C, V)$. In other words, it is an element in the abelian group
\[
\Hom_{\Gamma} (\Delta^0_C, V) := \{\varphi \in \Hom (\Delta^0_C, V) \:\: \vline \:\: \varphi (\gamma \cdot x) |\!|_{\gamma} = \varphi(x) \tu{ for all } \gamma \in \Gamma, x \in \Delta^0_C\}, 
\] where $|\!|_{\gamma}$ denotes the right action of $\gamma$ on $V$. We will denote this space by 
\[\Symb_{\Gamma, C} (V) := \Hom_{\Gamma} (\Delta^0_C, V).\]
\ed

\beg \label{example:the_set_C}
The primary examples of interest are the following.
\be
\item $C = \P^1 (\Q)$. This is the case of the \textit{$V$-valued (full) modular symbols}. As in the introduction, we write $\Delta^0$ for $\Delta^0_{\P^1(\Q)}$ and $\Symb_{\Gamma} (V)$ for $\Symb_{\Gamma, \P^1(\Q)} (V).$ 
\item Let $\Gamma \subset \Gamma_0 (N)$ be a subgroup where $N$ is divisible by $l^2$ for some prime number $l$. Then let $C = \P^1 (\Q) - \Gamma_0 (N) \cdot\{\infty\}$ be the complement of $\Gamma_0(N)$-translates of $\infty \in \P^1(\Q)$. One can check that $C$ is $\Gamma$-invariant. 
\item $C = \Gamma \{0\} \sqcup \Gamma \{\infty \}$, where $ \Gamma_1(p^m N) \subset \Gamma \subset \Gamma_0 (p^m N)$ for some $N \ge 1$ and $m \ge 2$. 
\ee 
\eeg 

\br[Hecke actions I]  
 
Let $\Sigma_0 (p) \subset \GL_2 (\Q)$ be the submonoid given by
\begin{equation*} 
\Sigma_0 (p)= \left\{ \begin{pmatrix}a & b \\ c & d \end{pmatrix} \in M_2 (\Z_{(p)}) \: \vline \: p \nmid a, p \mid c, ad - bc \ne 0 \right\}.
\end{equation*}
If the action of $\Gamma$ on $V$ extends to an action of $\Sigma_0 (p)$ and moreover if $C$ is stabilized by the action of $\Sigma_0 (p)$, then one may equip $\Symb_{\Gamma, C} (V)$ with the \emph{Hecke actions} induced from the action of $\Sigma_0 (p)$. 
\er 

\br[Hecke actions II] \label{remark:Hecke_action}
In practice, we will employ a slightly weaker form of Hecke actions. More precisely, suppose that we have
\begin{equation} \label{eq:form_of_Gamma}
\Gamma_1(p^m N) \subset \Gamma  \subset \Gamma_0 (p^m N),
\end{equation} where $m \ge 1$ and $(p, N) = 1$.  Let us assume that the action of $\Gamma$ on $V$ extends to an action of $\Sigma_0 (p)$, and further assume that $C$ is stabilized by the action of $ \big(\begin{smallmatrix} 1 & 0 \\ 0 & v \end{smallmatrix}\big)$ for all prime numbers $v$ satisfying $v \nmid N$, then we have induced  double coset actions $[\Gamma \big(\begin{smallmatrix} 1 & 0 \\ 0 & v \end{smallmatrix} \big) \Gamma ]$ for $v \nmid N$ on the modular symbols $\Symb_{\Gamma, C} (V)$. We again refer to these actions by \textit{Hecke actions} 
and denote them by $T_v = [\Gamma \big(\begin{smallmatrix} 1 & 0 \\ 0 & v \end{smallmatrix} \big) \Gamma ]$ for $v \nmid pN$ and $U_p =  [\Gamma \big(\begin{smallmatrix} 1 & 0 \\ 0 & p \end{smallmatrix} \big) \Gamma ]$, as usual. 
\er  

\br[Diamond operators] \label{remark:diamond}
In the setup of Remark \ref{remark:Hecke_action}, if moreover the set $C $ is preserved by the action of $\Gamma_0 (p^m N)$, 
then we also have well-defined diamond operators $\{\gr{d} | d \in (\Z/p^m N \Z)^\times \}$  on $\Symb_{\Gamma, C} (V)$.
\er 

\beg \label{example:Hecke_action} 
As an important special case of Remark \ref{remark:Hecke_action}, let 
$C = \P^1 (\Q) - \Gamma_0 (N) \{\infty\}$ as in Example \ref{example:the_set_C}. Let $V$ be an abelian group with a right $\Sigma_0(p)$ action as in Remark \ref{remark:Hecke_action} and let $\Gamma$ be as in (\ref{eq:form_of_Gamma}). Let us check that the assumptions in the remarks above are satisfied. First note that $C$ is preserved by $\Gamma_0 (p^m N)$. 
Next, for a prime number $v$, write $\gamma_{v} = \big(\begin{smallmatrix} 1 &  0 \\ 0 & v \end{smallmatrix} \big)$, then  $\gamma_v$ preserves $C$ for each $v \nmid N$. To see this, suppose that we have  $\gamma_{v} \cdot z = \gamma \cdot \infty$ for some $\gamma = \big(\begin{smallmatrix} a & b \\ cN & d \end{smallmatrix} \big) \in \Gamma_0 (N)$, then 
\[
z = \gamma_{v}^{-1} \gamma \cdot \infty = \frac{a v }{c N} \in \Gamma_0 (N) \cdot \{\infty\}
\]
since $(av, N) = 1$.  Therefore we have well-defined Hecke operators $U_p$, $T_v$ for $v \nmid pN$, and diamond operators $\gr{d}$ for $d \in (\Z/p^m N\Z)^\times$ on $\Symb_{\Gamma, C} (V)$.    
\eeg

\subsection{Relation with (classical) modular forms}\label{ss:Eichler_Shimura}
In this subsection, assume that $C$ is $\Gamma$-invariant and carries a Hecke action as in Remark \ref{remark:Hecke_action}. The short exact sequences 
\[ 
0 \ra \Delta^0_C \ra \Delta_C \xrightarrow{\deg} \Z \ra 0 
\]
for varying choices of $C$ induce the following commutative diagram with exact columns and rows (see \cite[Section 2.3]{BD}): 
\begin{equation} \label{eq:les_commutative_diagram} 
\begin{tikzcd}[column sep = 1em, row sep = 1.2em]
& \ker (\tu{res}_C) \arrow[r, equal] \arrow[d, hook] &   \ker (\tu{res}_C) \arrow[d, hook]  \\
0 \arrow[r] & \tu{BSymb}_{\Gamma} (V)  \arrow[r] \arrow[d, two heads] & \Symb_{\Gamma} (V) \arrow[r, "h"] \arrow[d, "\tu{res}_C"]  & H^1 (\Gamma, V) \arrow[r] \arrow[d, equal]  & \tu{coker} (h) \arrow[r] \arrow[d]  & 0 \\
0 \arrow[r] & \tu{BSymb}_{\Gamma, C} (V)  \arrow[r] & \Symb_{\Gamma, C} (V) \arrow[r, "h_C"] & H^1 (\Gamma, V) \arrow[r] & \tu{coker} (h_C) \arrow[r] & 0
\end{tikzcd}
\end{equation}
where  $\tu{Bymb}_{\Gamma, C}(V)$ denotes the \emph{boundary partial modular symbols}     \[\tu{BSym}_{\Gamma, C}(V) :=  \Hom_{\Gamma} (\Delta_C, V)/V^{\Gamma}\]
and $\tu{BSymb}_{\Gamma} (V)$ denotes $\tu{BSymb}_{\Gamma, \P^1(\Q)} (V)$.  
Let us fix a coefficient field $L$ of characteristic $0$ for the rest of this subsection. Let $k \ge 0$ be an integer and let $\tu{P}(L)^{\deg \le k}$ denote the space of polynomials of one variable $z$ of degree $\le k$ with coefficients in $L$, equipped with a left action of $\GL_2 (\Q)$ by\footnote{Here we adopt the convention of \cite{Pollack_Stevens} which is slightly different from that of \cite{BD}.} 
\begin{equation}\label{eq:action_on_P_L}
    (\gamma \cdot f )  (z) = (a + cz)^k f\Big(\frac{b + dz}{a + cz}\Big).
\end{equation}
Let $V^k = V^k (L)$ denote the linear dual of $P(L)^{\deg \le k}$, which is naturally acted upon by $\GL_2(\Q)$ on the right by $\psi |\!|_{\gamma} (f) = \psi (\gamma \cdot f)$ for all $f \in P(L)^{\deg \le k}$.

\bp[Eichler--Shimura] \label{prop:Eichler_Shimura}   
Suppose that we are in the setup of Remark \ref{remark:Hecke_action} and Remark  \ref{remark:diamond}. In particular, we have 
\[\Gamma_1(p^m N) \subset \Gamma\subset \Gamma_0(p^m N)\] for some $m \in \Z_{\ge 1}$ and $(p, N) = 1$. 
Consider the diagram (\ref{eq:les_commutative_diagram}) for $V = V^k(L)$. Let $M_{k+2}(\Gamma, L)$ (resp. $S_{k+2} (\Gamma, L)$, resp. $E_{k+2} (\Gamma, L)$) denote the space of modular forms (resp.  cusp forms, resp. Eisenstein series)  of weight $k+2$ and level $\Gamma$.
Then after possibly replacing $L$ by a finite extension, the following holds:   
\be
\item There are isomorphisms 
\begin{equation} \label{eq:hecke_equivariant_ES}
 \im (h) \isom S_{k+2}(\Gamma, L)^{\oplus 2}  \quad \tu{ and} \quad   \tu{BSymb}_{\Gamma} (V^k) \isom E_{k+2} (\Gamma, L),
\end{equation} 
which are Hecke equivariant with respect to the action of $U_p$, $T_v$ for $v \nmid pN$, and the diamond operators $\gr{d}$ for $d \in (\Z/p^m N\Z)^\times$.  
\item  Let $\T$ denote the Hecke algebra generated by $U_p$, $T_v$ for $v \nmid pN$, and the diamond operators $\gr{d}$ for $d \in (\Z/p^m N\Z)^\times$. Then there are natural $\T$-equivariant inclusions 
\[
S_{k+2} (\Gamma, L)^{\oplus 2} \subset   (\Symb_{\Gamma, C} (V^k))^{\tu{s.s.}}
\subset M_{k+2}(\Gamma, L)^{\oplus 2} 
\]
where the middle term denotes the the semi-simplification of $\Symb_{\Gamma, C} (V^k)$ with respect to the action of  $\T$. 
\ee 
\ep

\bproof 
The first claim is the Eichler--Shimura isomorphism (see \cite[Proposition 2.5]{BD}). The second claim follows from the first by a simple diagram chase (note that in particular we have $\im (h) \subset \im (h_C) \subset   (\Symb_{\Gamma, C} (V^k))^{\tu{s.s.}}$, see also \cite[Corollary 2.7]{BD}). 
\eproof

\subsection{$\Gamma$-structure of $\Delta_C^0$} \label{ss:Delta_free_over_Gamma} In this subsection we prove some important structural results on $\Delta_C^0$. Let us first recall some results of \cite{Pollack_Stevens} in the case where $C = \P^1(\Q)$. Let $\mH$ denote the upper half plane and write $\mH^* = \mH \cup \P^1(\Q)$ with its usual topology. Let $\Gamma \subset \PSL_2(\Z)$ be a finite index subgroup without $3$-torsion. Recall that a \emph{modular chain} is a finite sum of geodesic paths in $\mH^*$ of the form $e_{\gamma}$ where $\gamma \in \GL_2(\Z)$ (see introduction for the notation). Note that $\Delta^0$ can be identified with the relative homology $H^1 (\mH^*, \P^1(\Q); \Z)$, which may be regarded as the quotient of the abelian group of modular chains with endpoints on $\partial \mH^* = \P^1(\Q)$ by the boundary modular chains\footnote{Here a boundary modular chain means a modular chain that forms a closed loop (a $1$-cycle) in $\mH^*$}. The continued fractions trick implies that there is a $\PSL_2(\Z)$-equivariant surjection \[
\pi: \Z[\PSL_2(\Z)] \twoheadrightarrow \Delta^0
\]
sending $[\big(\begin{smallmatrix} a & b \\ c & d \end{smallmatrix} \big)] \mapsto \{\frac{a}{b}\} - \{\frac{c}{d}\} \in \Delta^0$, where the latter may be viewed as the (geodesic) path from $\frac{c}{d}$ to $\frac{a}{b}$ in $H^1(\mH^*,\P^1(\Q))$ and called a \emph{unimodular path}. Manin showed in \cite{Manin} that $\ker (\pi)$ is the left ideal of $\Z[\PSL_2(\Z)]$ generated by two elements $[1] + [\tau] + [\tau^2]$ and $[1] + [\sigma]$ where $\tau = \big(\begin{smallmatrix} 0 & -1 \\ 1 & -1 \end{smallmatrix} \big) $ and $\sigma = \big(\begin{smallmatrix} 0 & -1 \\ 1 & 0 \end{smallmatrix} \big)$. 

\bp[Pollack--Stevens] \label{prop:PS_shape_of_domain} Let $\Gamma \subset \PSL_2(\Z)$ be a finite index subgroup without $3$-torsion. There exists a fundamental domain $\mF$ for the action of $\Gamma$ on $\mH^*$, such that $\mF$ is contained in the strip $\{0 \le \tu{Re}(z) \le 1\}$, the vertices of $\mF$ are all cusps, and the boundary $\partial \mF$ of $\mF$ is a union of unimodular paths.
Moreover, let $\mE$ be the set of all oriented unimodular paths passing through the interior of $\mF$ together with the ones contained in $\partial \mF$ (where $\partial \mF$ is oriented clockwise). Then any oriented unimodular path is $\Gamma$-equivariant to a unique element in $\mE$. 
\ep  

\bproof 
This is \cite[Theorem 2.3 and Lemma 2.5]{Pollack_Stevens}. 
\eproof 

\beg \label{example:fund_domain} 
Let $\mT$ be the ideal triangle with vertices $(0, \infty, \frac{1+\sqrt{-3}}{2})$ and let 
\[\mR = \mT \cup \tau \mT \cup \tau^2 \mT\] where $\tau = \big(\begin{smallmatrix} 0 & -1 \\ 1 & -1 \end{smallmatrix} \big) $  as before. Namely, $\mR$ is the ideal triangle with vertices $(0, \infty, 1)$. More generally, write $\mR_{(a, b, c)}$ for the ideal triangle with vertices $(a, b, c)$. 
\be
\item Consider 
\[\PSL_2 (\Z) = \Big(\bigsqcup_{i = 0}^{4} \Gamma_0(5) \gamma_i \Big)\sqcup \Gamma_0 (5) \gamma_\infty\] where $\gamma_i =  \big(\begin{smallmatrix} 1 & 0 \\ i & 1 \end{smallmatrix} \big)$ for $i = 0, ..., 4,$ and $\gamma_\infty =  \big(\begin{smallmatrix} 0 & -1 \\ 1 & 0 \end{smallmatrix} \big)$. Observe that $\Gamma_0 (5) \tau = \Gamma_0 (5) \gamma_4$ and that $\Gamma_0 (5) \tau^2 = \Gamma_0 (5) \gamma_{\infty}$. A further computation shows that we can write 
\[\PSL_2 (\Z) = \bigsqcup_{j = 0}^2 \Big(\Gamma_0(5) \sqcup \Gamma_0 (5) \gamma_1 \Big) \cdot \tau^j.\]
Therefore,   the union $\mR \cup \gamma_1 \mR$ is a  fundamental domain for $\Gamma_0(5)$, where $\gamma_1 \mR = \mR_{(0, 1, \frac{1}{2})}$. 
\item A similar computation shows that one fundamental domain for $\Gamma_0 (11)$ can be chosen to be the union $\mR \cup \mR_{(0, 1, \frac{1}{2})} \cup \mR_{(0, \frac{1}{2}, \frac{1}{3})} \cup \mR_{(\frac{1}{2}, 1, \frac{2}{3})}$. 
\ee 
\eeg

Let us further assume that $\Gamma$ contains no $2$-torsion elements. Now we fix a choice of a fundamental domain $\mF$ for $\Gamma$ as in Proposition \ref{prop:PS_shape_of_domain}.  For each unimodular path $e \in \partial \mF$, let $\cl e$ denote the path obtained from $e$ by switching the orientation, so that we have $e + \cl e = 0$ in $\Delta^0$. By proposition \ref{prop:PS_shape_of_domain}, $\cl e$ is $\Gamma$-equivariant to a unique unimodular path $e^* \in \partial \mF$ (for example, we have $\cl e_\infty = \big(\begin{smallmatrix} 1 & 1 \\ 0 & 1 \end{smallmatrix} \big) e_\infty^*$, where $e_\infty  = \{1\} -\{\infty\}$ and $e_\infty^* = \{\infty\} - \{0\}$). Moreover, we know that $e \ne e^*$ since $\Gamma$ contains no $2$-torsion.   
Therefore, reversing orientation decomposes $\partial \mF$ into orbits of size $2$. Let us enumerate the elements in $\partial \mF$  as 
\begin{equation} \label{eq:generators_for_Delta}
e_\infty, e_{\infty}^*,  e_1, e_1^*, \dots, e_{t}, e_t^*  
\end{equation}
for some integer $t$, and let $\gamma_i \in \Gamma$ be the element such that $\cl e_i = \gamma_i e_i^*$ for each $i = 1, ..., t, \infty$.

\bp[Pollack--Stevens]\label{prop:PS_generator_and_relations}
Suppose that $\Gamma$ is a finite index subgroup of $\PSL_2(\Z)$ without any $2$- or $3$-torsion elements. Let $e_i$ and $e_i^*$ be as above. Then $\Delta^0$ is generated as a $\Z[\Gamma]$-module by $e_1, ..., e_t, e_\infty$ with relations generated by a single equation 
\[
\sum_{i = \infty, 1, ..., t}(1 - \gamma_i^{-1}) e_i = 0. 
\]
\ep

\bproof 
This is (a special case of) \cite[Theorem 2.6]{Pollack_Stevens}. 
\eproof 

Note that the single relation satisfied by these generators precisely come from the loop $\partial \mF$, so one might hope that ``breaking the loop'' would yield a finite free $\Z[\Gamma]$-submodule of $\Delta^0$. Let us formulate this as the following: 

\begin{question} \label{question:Delta^0_C_structure}
Keep the assumption on $\Gamma$ from Proposition \ref{prop:PS_generator_and_relations}. Let $C \subsetneq \P^1(\Q)$ be a $\Gamma$-invariant proper subset. When is $\Delta^0_C$ a free module over $\Z[\Gamma]$?  
\end{question}

In the rest of this subsection, we provide an answer to Question \ref{question:Delta^0_C_structure} in some special cases that would suffice for our applications to the Halo conjecture. 

\bl \label{lemma:unique_translate_of_infinity}
Let $l$ be a prime number satisfying $l \equiv 11 \!\mod 12$. Let $\mF$ be a fundamental domain for $\Gamma_0 (l^2)$ as in Proposition \ref{prop:PS_shape_of_domain} and let $V(\mF) \subset \P^1 (\Q)$ be the set of vertices of $\mF$. Then 
\[V(\mF) \cap  \big(\Gamma_0 (l^2)\cdot\{\infty\}\big) = \{\infty\}\] 
where the intersection is taken inside $\P^1(\Q)$. 
\el 

\bproof The assumption on $l$ implies that  $\Gamma_0(l)$ contains no $2$- or $3$-torsion elements. 
As in Example \ref{example:fund_domain}, we write\footnote{The middle term in (\ref{eq:expression_of_fd_for_l}) is introduced to simply notation later on in the proof.} 
\begin{equation} \label{eq:expression_of_fd_for_l}
\PSL_2 (\Z) = \Gamma_0 (l) \Big(\bigsqcup_{i = \infty, 0, ..., l-1}\gamma_i   \Big) := \Big(\bigsqcup_{i = 0}^{l-1} \Gamma_0(l) \gamma_i \Big)\sqcup \Gamma_0 (l) \gamma_\infty 
\end{equation} where $\gamma_i =  \big(\begin{smallmatrix} 1 & 0 \\ i & 1 \end{smallmatrix} \big)$ for $i = 0, ..., l-1,$ and  $\gamma_\infty =  \big(\begin{smallmatrix} 0 & -1 \\ 1 & 0 \end{smallmatrix} \big)$. Note that we have \[\Gamma_0 (l) \tau = \Gamma_0 (l) \gamma_{l-1} \quad \text{ and} \quad \Gamma_0 (l) \tau^2 = \Gamma_0 (l) \gamma_{\infty}.     
\]
Therefore, the decomposition in (\ref{lemma:unique_translate_of_infinity}) can be rewritten as 
\begin{equation*}
  \label{eq:expression_of_fd_for_l}
\PSL_2 (\Z) = \Gamma_0 (l) \Big(\bigsqcup_{j=0}^{(l-2)/3} \gamma_{r_j}   \Big) \Big( \tu{id} \sqcup \tau \sqcup \tau^2   \Big)  
\end{equation*}
where $r_j$ are (distinct) integers among $\{0, ..., l-2\}$ with $r_0 = 0$ (so $\gamma_{r_0} = \tu{id}$). This further implies that 
\begin{equation}
  \label{eq:expression_of_fd_for_l_square}
\PSL_2 (\Z) = \Gamma_0 (l^2)  \Big(\bigsqcup_{k=0}^{l-1} \beta_k   \Big)  \Big(\bigsqcup_{j=0}^{(l-2)/3} \gamma_{r_j}   \Big) \Big( \tu{id} \sqcup \tau \sqcup \tau^2   \Big)  
\end{equation}
where  $\beta_k = \gamma_{lk} =  \big(\begin{smallmatrix} 1 & 0 \\ lk & 1 \end{smallmatrix} \big)$ for $k = 0, ..., l-1$. Therefore, the fundamental domain is the union of certain $\Gamma_0(l^2)$-translates of the ideal triangles $\beta_k \gamma_{r_j} \cdot \mR$ where $\mR = (\tu{id}\sqcup \tau \sqcup \tau^2) \mT = \mR_{(0, \infty, 1)}$ from Example \ref{example:fund_domain}. Note that for $k=j=0$, the corresponding $\Gamma_0(l^2)$-translate is precisely $\mR$. Thus it suffices to show that, for all $k$ and $j$ such that $k^2 + r_j^2 \ne 0$, the $\beta_k \gamma_{r_j}$-translates of the vertices $(0, 1, \infty)$ do not live in $\Gamma_0 (l^2) \{\infty\}$. But this is clear, since we have 
\[
\beta_k \gamma_{r_j} \cdot 0 = 0, \quad \beta_k \gamma_{r_j} \cdot \infty = \frac{1}{lk + r_j}, \quad \beta_k \gamma_{r_j} \cdot 1 = \frac{1}{lk + r_j + 1},
\]
none of which belong to $\Gamma_0(l^2) \cdot\{\infty\}$ (for the middle term we use that $k$ and $r_j$ cannot both be $0$ and for the last term we use that $r_j \ne l-1$). 
\eproof 

\bp  \label{prop:partial_free_for_l_square}
Retain assumptions from Lemma \ref{lemma:unique_translate_of_infinity} and let $C = \P^1(\Q) - \Gamma_0 (l^2) \cdot \{\infty\}$. Let $\{e_\infty, e_\infty^*, e_1, e_1^*, ..., e_t, e_t^*\}$ be the oriented unimodular paths in $\partial \mF$ as in (\ref{eq:generators_for_Delta}). Then the elements $e_1, ..., e_t$ provide a set of free generators of $\Delta_C^0$ as a $\Z[\Gamma_0 (l^2)]$-module. 
\ep 

\bproof 
By Lemma \ref{lemma:unique_translate_of_infinity} we know that $e_1, ..., e_t$ belong to $\Delta^0_C$. Moreover, by Proposition \ref{prop:PS_generator_and_relations} there is no relation satisfied among these elements, thus it suffices to show that they generate $\Delta^0_C$. For this, let $z = \{a_1\} - \{a_0\}$ be a modular path in $\Delta^0_C$, and we want to show that $z$ can be written as a sum of $\Gamma_0 (l^2)$-translates of the paths $e_i$'s. Using the continued fractions trick, we may write $z$ as a sum of unimodular paths 
\begin{equation} \label{eq:cont_frac_trick}
z = (\{a_1\} - \{a_2\}) + (\{a_2\} - \{a_3\}) + \dots + (\{a_s\} - \{a_0\}).  
\end{equation}
Write $p_i = \{a_i\} - \{a_{i+1}\}$ for the $i^{th}$-unimodular path appearing in (\ref{eq:cont_frac_trick}) (where $a_{s+1}$ is set to be $a_0$), then each $p_i$ is $\Gamma_0(l^2)$-equivaraint to a unique element 
$\sq p_i \in \mE$ by Proposition \ref{prop:PS_shape_of_domain}. 
Let $p_{k_1}$ be the first path among the $p_i$'s such that $\sq p_{k_1}$ goes through the point $\infty$ (note that if such $p_{k_1}$ does not exist then each $\sq p_i$ is a sum of paths from the set $\{e_1, e_1^*, ..., e_t, e_t^*\}$ and we are done). First observe that we must have $a_{k_1 + 1} \in \Gamma_0(l^2) \cdot \{\infty\} $ (in other words, $a_{k_1} \notin \Gamma_0 (l^2) \cdot \{\infty\}$), or equivalently, $\sq p_{k_1}$ is of the form $\{b_{1}\} - \{\infty\}$.  This follows from the minimality of the index $k_1$. More precisely, suppose that instead we have 
$\sq p_{k_1} = \{\infty\} - \{b'\}$ for some $b' \in \Q$, then, as $a_1 \notin \Gamma_0 (l^2) \{\infty\}$ we know that $k_1 \ne 1$, thus $\sq p_{k_1 -1}$ is of the form $\{c\} - \{c'\}$ where $c' \in \Gamma_0 (l^2) \{\infty\}$. This forces $c' = \infty$ by Lemma \ref{lemma:unique_translate_of_infinity}, and thus contradicts the minimality of $k_1$. Therefore, indeed we have 
\begin{equation} \label{eq:tilde_p_k_1}
    \sq p_{k_1} = \delta \cdot p_{k_1}= \{\delta \cdot a_{k_1} \} - \{\delta \cdot a_{k_1+1}\} =  \{b_{1}\} - \{\infty\}
\end{equation}
for some $\delta \in \Gamma_0(l^2)$ and $b_1 \in \Q$. Note that $k_1 \ne s$ since $a_0 \notin \Gamma_0 (l^2) \{\infty\}$. Hence, $2 \le k_1 \le s-1$ and we have 
\begin{equation} \label{eq:tilde_p_k_2}
    \sq p_{k_1+1} = \delta' \cdot p_{k_1 + 1} =  \{\delta' \cdot a_{k_1+1} \} - \{\delta' \cdot a_{k_1+2}\} = \{\infty\} - \{c_1\} \in \mE.
\end{equation}
Now applying Lemma \ref{lemma:unique_translate_of_infinity} one more time we know that $\delta' \cdot a_{k+1} = \infty$. Now note that we have $\delta' \delta^{-1}: \infty \mapsto \infty$, which implies that $\delta' \delta^{-1} =   \big(\begin{smallmatrix} 1 & k \\ 0 & 1 \end{smallmatrix} \big) $ for some $k \in \Z$, in other words, we have $\delta' =  \big(\begin{smallmatrix} 1 & k \\ 0 & 1 \end{smallmatrix} \big) \cdot \delta  \in \PSL_2(\Z)$. Combining (\ref{eq:tilde_p_k_1}) and (\ref{eq:tilde_p_k_2}), we have 
\begin{align}
 p_{k_1} + p_{k_1 +1 } & = \delta^{-1} \cdot (\sq p_{k_1} + \sq p_{k_1+1}) + \delta^{-1} \cdot \big(  \big(\begin{smallmatrix} 1 & -k \\ 0 & 1 \end{smallmatrix} \big) \cdot \sq p_{k_1+1} - \sq p_{k_1+1} \big) \\
 \nonumber & =  \delta^{-1} \cdot (\{b_1\} - \{c_1\})  + \delta^{-1} \cdot \big( \{c_1\} - \{c_1 - k\}  \big), 
\end{align}
which clearly can be generated by  $\Gamma_{0}(l^2)$-translates of the $e_i$'s for $i = 1,..., t$. Repeating this procedure, we obtain the desired proposition. 
\eproof 

\bc \label{cor:free_over_Gamma} Retain assumptions from Proposition \ref{prop:partial_free_for_l_square} and suppose that $l \ne p$. Suppose that 
$ \Gamma \subset \Gamma_0 (l^2 p^m)$ for some $m\geq 1$. Then the choice of generators $e_1, \ldots, e_t$ in Proposition \ref{prop:partial_free_for_l_square} and  coset representatives of $\Gamma_0 (l^2)/\Gamma$ induces an isomorphism 
\[
\Symb_{\Gamma, C} (V) \cong V^{\oplus st}
\]
where $s$ is the index of $\Gamma \subset \Gamma_0 (l^2)$. 
\ec 

\bproof 
This is straight-forward. 
\eproof

\section{Overconvergent (partial) modular symbols} \label{section:oc_symb}
In this section we discuss overconvergent partial modular symbols and an ``integral model'' for them. 
\subsection{Overconvergent partial modular symbols}  \label{ss:oc_partial_symbols} 
We largely follow the conventions from \cite{Pollack_Stevens}. For each  $r \in |\C_p^\times|$, we consider the rigid analytic space 
\[B_{r} = \{z \in {\C_p} \: \vline \:  |z - a| \le r \tu{ for some } a \in \Z_p\}.\] For a $\Q_p$-Banach algebra $L$, let $A_r(L)$ be the $L$-Banach module of functions $f: \Z_p \ra L$ that are locally analytic of radius $r$ (namely, the restriction to $\Z_p$ of an analytic function on $B_r$), equipped with the supremum norm. Let $D_r (\Q_p)$ be the (strong) Banach dual of $A_r(\Q_p)$ and let $D_r(L) = D_r(\Q_p) \widehat \otimes_{\Q_p} L$. We further define 
\[\mA^{\dagger}_r (L) := \varinjlim_{s > r} A_s (L), \qquad \text{and}\quad  \mD^{\dagger}_r(L) : = \varprojlim_{s > r} 
D_s (L) \: \: \:  \] 
with the inductive (resp. projective) limit topology. We also write 
\[\mA^{\dagger}(L) := \varinjlim_{s > 0} A_s (L),    \qquad \textup{and} \quad \mD^{\dagger}(L)  := \varprojlim_{s > 0} 
D_s (L).\footnote{In other words, $\mA^{\dagger}(L)= \mA^{\dagger}_0(L)$ and $\mD^{\dagger}(L) = \mD^{\dagger}_0(L).$ This notation slightly differs  from \cite{Pollack_Stevens} but agrees with \cite{BD}.}\] 
When the context is clear, we omit ``$L$'' and just write $A_r, D_r, \mA^{\dagger}_r, \mD^{\dagger}_r$, etc.
Note that, we have natural inclusions
\begin{equation} \label{eq:inclusions}
\mD^{\dagger} \hookrightarrow D_r \hookrightarrow \mD_{r}^{\dagger}  \hookrightarrow D_s 
\end{equation}
for $0 < r < s$ (see  \cite[Section 3.1]{Pollack_Stevens_critical}). For a continuous character $\chi: \Z_p^\times \ra L^\times$ which is locally analytic of radius $r_0$, then for any $r \le r_0$ we consider the action of $\Sigma_0 (p)$ on $A_{r}$ defined as follows: for $\gamma =  \big(\begin{smallmatrix} a & b \\ c & d \end{smallmatrix} \big)  \in \Sigma_0 (p)$, we let
\begin{equation} \label{eq:action_on_A_r}
    (\gamma \cdot f) (z) = \chi(a+cz) f(\frac{b+dz}{a+cz}),
\end{equation}
where $\chi(a+cz)$ denotes the locally analytic map $\Z_p \ra L$ given by $z \mapsto \chi(a+cz)$. 
We then obtain a right action on $D_r$ by 
\begin{equation} \label{eq:action_on_D_r}
   \psi |\!|_\gamma  (f) = \psi (\gamma \cdot f). 
\end{equation}
This naturally induces actions on $\mA^{\dagger}_r$ and $\mD^{\dagger}_r$ for all $r < r_0$, as well as $ \mA^{\dagger}$ and $\mD^{\dagger}$. We write $A_r^\chi$ (resp. $D_r^{\chi}$, resp. $\mA^{\dagger, \chi}_r$, resp. $\mD^{\dagger, \chi}_r$, resp. $\mA^{\dagger, \chi}$, resp. $\mD^{\dagger, \chi}$) for the corresponding $L$-Banach (or -Frechet) module equipped with the ``$\chi$-twisted'' action of $\Sigma_0 (p)$ defined above.

\br  \label{remark:transpose} For psychological comfort, we remark that the right action of $\Sigma_0 (p)$ on $D_r^{\chi}$ and its variants obtained from the left action on $A_r^{\chi}$ 
can be thought of as a ``transpose action" in the following sense: \emph{if}  $\{f_1, f_2, ...\}$ form an $L$-basis of $A_r$, write $f_i^{\vee} \in D_r$ for the dual linear functionals. 
If the action of $\gamma$ on $A_r$ is  described by an infinite matrix $P_{m,n}$, in other words, $\gamma \cdot f_n = \sum_m P_{m,n} f_m$. Then the action of $\gamma$ on the $L$-span of the linear functionals $\{f_i^\vee\}$'s is given by 
\[
f_j^{\vee}|\!|_\gamma (f_i) = f_j^{\vee} (\sum_{m} P_{m i} \; f_m)  = f_j^{\vee} (P_{j i} \; f_j) = P_{ji}.
\]
In other words, $f_j^{\vee}|\!|_\gamma = \sum_{i} P_{ji} \; f_i^{\vee}$. Thus ``the action of $\gamma$ on $\{f_i^\vee\}$'s is obtained by  multiplication by the transpose matrix $P_{m, n}^{t}$ on the right.''
\er 

\begin{assumption} \label{assumption:sec3}
For the rest of Section \ref{section:oc_symb} we assume that $C$ is a $\Gamma_0 (Np)$-invariant subset of $\P^1(\Q)$ where $p \nmid N$. Let us further assume that $C$ is stabilized by the action of $ \big(\begin{smallmatrix} 1 & 0 \\ 0 & v \end{smallmatrix}\big)$ for all primes $v$ satisfying $v \nmid N$.  
\end{assumption}
In particular, for a subgroup $\Gamma \subset \Gamma_0 (pN)$ as specified in (\ref{eq:form_of_Gamma}) and an abelian group $V$ equipped with a right action by $\Sigma_0(p)$, we have well-defined Hecke operators $U_p,$ $T_v$ for all $v \nmid pN$, and $\gr{d}$ for $d \in (\Z/p^m N\Z)^\times$ on $\Symb_{\Gamma, C} (V)$ by Remark \ref{remark:Hecke_action} and Remark  \ref{remark:diamond}. 
\bd 
Let $\chi: \Z_p^\times \ra L^\times$ be a continuous character locally analytic of radius $r_0$ as above. \emph{The space of overconvergent partial modular symbols supported on $C$ of weight $\chi$ and tame level $\Gamma_0 (N)$} (resp. \emph{and locally analytic of radius} $r$ where $r < r_0$) is  the space 
\[ \Symb_{\Gamma_0(Np), C} (\mD^{\dagger, \chi})
 \quad (\tu{resp. }  \Symb_{\Gamma_0 (Np), C} (\mD^{\dagger, \chi}_{r}) \:  ). 
\]
\ed  
For the next lemma, suppose that $L$ is a field extension of $\Q_p$. 
Recall that a continuous character $\chi: \Z_p^\times \ra L^\times$ is called a \emph{classical weight character} if it is of the form $\chi = (k, \epsilon)$ given by 
\[\chi(z) = {(k, \epsilon)} (z) := z^k \epsilon (z),\] 
where $k \in \Z$ and $\epsilon$ is a finite character $\Z_p^\times \ra (\Z/p^m \Z)^\times \ra L^\times$. We call $\epsilon$ the $p$-\emph{Nebentypus character} of $(k, \epsilon)$, and the smallest such $m$ is called the \emph{conductor} of $(k, \epsilon)$. Next observe that, for each $k \ge 0$ and $r \in |\C_p^\times|$, the natural inclusion $P(L)^{\deg \le k} \hookrightarrow \mA^{\dagger}_r$ induces a surjective map $ \mD_r^{\dagger} \twoheadrightarrow V^k = V^k(L)$ (The spaces $P(L)^{\deg \le k}$ and $V^k$ are defined in Subsection \ref{ss:Eichler_Shimura}). 

\bl \label{lemma:epsilon_component} Let ${(k, \epsilon)}$ be  a classical weight character of conductor $m$. Then for $r < \frac{1}{p^m}$,  the map $ \mD_r^{\dagger} \twoheadrightarrow V^k$ above induces a natural \emph{specialization} map \begin{equation} \label{eq:specialize_to_classical}
    \Symb_{\Gamma_0 (Np), C} (\mD_r^{\dagger, (k, \epsilon)}) \lra \Symb_{\Gamma_0(N) \cap \Gamma_1(p^m), C} (V^k)_{\epsilon},
\end{equation}
equivariant for the Hecke operators $U_p$ and $T_v$ where $v \nmid Np$. Here the subscript on the right hand side denotes the $\epsilon$-part of the $\Gamma_0(p^m)$-action; in other words, it consists of modular symbols $\phi \in \Symb_{\Gamma_0(N) \cap \Gamma_1(p^m), C} (V^k)$ satisfying $\gamma \cdot \phi = \epsilon(d) \phi$ for all $\gamma=\big(\begin{smallmatrix} a & b \\ c & d \end{smallmatrix} \big) \in \Gamma_0 (Np^m)$. These specialization maps are compatible when $r$ varies. 
\el 

\bproof 
For an element $f \in P(L)^{\le \deg k}$ and $\gamma = \big(\begin{smallmatrix} a & b \\ c & d \end{smallmatrix} \big) \in \Gamma_0(p^m)$, we will  use $\gamma\cdot^{k} f$ to denote the action of $\gamma$ on $f$ defined in (\ref{eq:action_on_P_L}), to distinguish it from the ``twisted action" defined by 
\begin{equation} \label{eq:action_twisted_by_epsilon}
\gamma\cdot^{(k,\epsilon)} f := \epsilon(a+cz)\cdot (\gamma\cdot^{k} f) = \epsilon(a+cz) (a + cz)^k f\Big(\frac{b + dz}{a + cz} \Big).
\end{equation}
By our assumptions on $\gamma$ and the conductor of $\epsilon$, the function $\epsilon(a+cz)$ is simply the constant function $\epsilon(a)$ in this case. Further note that $\epsilon(a+cz)  = 1$ for $\big(\begin{smallmatrix} a & b \\ c & d \end{smallmatrix} \big) \in \Gamma_1(p^m)$, so the map $\mD_r^{\dagger, (k, \epsilon)} \ra V^k$ is equivariant for the $\Gamma_1(p^m)$-action. Thus we have a natural map 
\[
 \Symb_{\Gamma_0 (Np), C} (\mD_r^{\dagger, (k, \epsilon)}) \lra \Symb_{\Gamma_0(N) \cap \Gamma_1(p^m), C} (V^k).
\]
Denote by $\cl \varphi$ the image of the $\varphi$ under this map, then for any  $\big(\begin{smallmatrix} a & b \\ c & d \end{smallmatrix} \big) \in \Gamma_0(p^m)$, we have 
\[
\gamma \cdot \cl \varphi (x) (f) = \cl \varphi(\gamma \cdot x) (\gamma\cdot^{k} f) = \epsilon(a)^{-1}  \varphi(\gamma \cdot x) (\gamma\cdot^{(k, \epsilon)} f)= \epsilon(d) \cl \varphi (x) (f) 
\]
for all $x \in \Delta_C^0$ and all $f \in P(L)^{\deg \le k}$. Here the last equality uses the fact that $\cl \varphi$ comes from $\varphi: \Delta^0_C \ra \mD_r^{\dagger, (k, \epsilon)}$ and is invariant under the action of $\Gamma_0(Np^m)$. This proves the lemma. 
\eproof 

\br \label{remark:epsilon_component} If we further assume that $N = l^2$ and $C = \P^1(\Q) - \Gamma_0 (N) \cdot \{\infty\}$ as in Proposition \ref{prop:partial_free_for_l_square}, then our proof in fact shows that the map (\ref{eq:specialize_to_classical}) factors as \begin{multline*}
\quad   \Symb_{\Gamma_0 (Np), C} (\mD_r^{\dagger, (k, \epsilon)}) \lra 
\Symb_{\Gamma_0 (Np^m), C} (\mD_r^{\dagger, (k, \epsilon)}) \\ 
\surjlra \Symb_{\Gamma_0 (Np^m), C} (V^{(k, \epsilon)}) \isom  \Symb_{\Gamma_0(N) \cap \Gamma_1(p^m), C} (V^k)_{\epsilon} \quad
\end{multline*}
where $V^{(k, \epsilon)}$ denotes the vector space $V^k$ but with a twisted right action by $\Sigma_0 (p)$ induced from the action $\gamma \cdot^{(k,\epsilon)} f$  in (\ref{eq:action_twisted_by_epsilon}). 
For the surjection in the middle, we use Corollary \ref{cor:free_over_Gamma}. 
\er 

Next let us record a lemma which essentially says that the overconvergent modular symbols descend to $\Gamma_0(p)$-level (while shrinking the radius of convergence at the same time). 
\bl \label{lemma:level_raising_and_lowering_at_p} Let $\chi$ be a continuous character locally analytic of radius $r_0$. Then for $r \le r_0$, there is a natural isomorphism  
\[
\Symb_{\Gamma_0 (Np^m), C} (\mD_{r/p}^{\chi}) \isom 
\Symb_{\Gamma_0 (Np^{m+1}), C} (\mD_r^{\chi}), 
\]
equivariant for the Hecke operators $U_p$ and $T_v$ for $v \nmid Np$. For $r < r_0$, the analogous statements hold for  $\mD_{r/p}^{\dagger, \chi}$ and  $\mD_{r}^{\dagger, \chi}$. 
\el 

\bproof 
This is proved by a similar argument as \cite[Lemma 4.4(4)]{Buzzard_auto_family}, which we recall in our context. Write $p B_r \subset B_{r/p}$ for the rigid space $\{z | z/p \in B_r \}$ and let $\mD^{\chi}_{p B_r}$ denote the Banach dual of the space of convergent functions on $p B_r$, equipped with the action of $\Sigma_0(p)$ given by the formula (\ref{eq:action_on_A_r}). 
Let $\Gamma^0 (p) \equiv \big(\begin{smallmatrix} * & 0 \\ * & * \end{smallmatrix} \big) \mod p$ consist of matrices that are lower triangular mod $p$.  Then we claim that there are natural Hecke equivariant isomorphisms
\begin{align*}
\Symb_{\Gamma_0 (Np^m), C} (\mD_{r/p}^{\chi}) &\isom  \Symb_{\Gamma_0 (Np^m)\cap \Gamma^0(p), C} (\mD_{p B_r}^{\chi})\\ & \isom \Hom_{\Gamma_0 (Np^{m+1})} (\Delta_C^0, \mD_{p B_r}^{\chi}|\!|_{\big(\begin{smallmatrix} 1 & 0 \\ 0 & p \end{smallmatrix} \big)}) \\
 & \isom  \Hom_{\Gamma_0 (Np^{m+1})} (\Delta_C^0, \mD_{r}^{\chi}).
\end{align*}
For the first arrow, both sides can be identified with the subspace of 
$\Symb_{\Gamma_0 (Np^m) \cap \Gamma^0 (p), C} (\mD_{r/p}^{\chi})$ that is invariant under the action of the matrix $\big(\begin{smallmatrix} 1 & 1 \\ 0 & 1 \end{smallmatrix} \big)$. The second isomorphism is given by sending $\varphi \longmapsto  \big(\begin{smallmatrix} 1 & 0 \\ 0 & p \end{smallmatrix} \big) \cdot \varphi$. The last isomorphism is straightforward. 
\eproof 

\subsection{Comparison theorem}\label{subsection: comparison theorem}
Keep Assumption \ref{assumption:sec3} from the previous subsection. Throughout Section \ref{subsection: comparison theorem}, let $L$ be a field  extension of $\Q_p$ and let $(k, \epsilon)$ be a classical weight of conductor $m$. 
\begin{notation}  Let $\Gamma \subset \Gamma_0 (Np)$ be a (congruence) subgroup. 
For any $h \in \R$, we use $\Symb_{\Gamma, C} (\mD_r^{\chi})^{<h}$ to denote the subspace of the Banach space $\Symb_{\Gamma, C} (\mD_r^{\chi})$ on which the $U_p$ operator acts with slope $< h$ (see  \cite[Definition 4.7]{Pollack_Stevens_critical}). For each $r$, let  
\begin{align*}
  \Symb_{\Gamma, C} (\mD_r^{\dagger, \chi})^{<h} &:= \Symb_{\Gamma, C} (\mD_r^{\dagger, \chi}) \cap \Symb_{\Gamma, C} (\mD_s^{\chi})^{<h}
\end{align*}
for some $s>r$, where the intersection takes place in  $\Symb_{\Gamma, C} (\mD_s^{\chi})$ (using (\ref{eq:inclusions})). We further write \[\Symb_{\Gamma, C} (\mD^{\dagger, \chi})^{<h} :=  \Symb_{\Gamma, C} (\mD^{\dagger, \chi}) \cap \Symb_{\Gamma, C} (\mD_s^{\chi})^{<h}  
\]
for some $s>0$. 
Note that these definitions \textit{a priori} depend on the choice of $s$. This dependence is removed by the next proposition.\footnote{Technically, Proposition \ref{prop:finite_slope_witness_no_p_level} only shows that the definitions above are independent of $s$ at classical weights. More generally, this dependence can be removed by a similar argument. For simplicity, we stick with our current setup, which suffices for the proof of the Halo conjecture.} 
\end{notation}

\bp[Pollack--Stevens] \label{prop:finite_slope_witness_no_p_level}
 The Banach space $\Symb_{\Gamma_0(Np), C}(\mD_r^{(k, \epsilon)})$ is orthonormalizable and the $U_p$-action is compact. Moreover, for any $h \in \R$ and any $s > r > 0$, the restriction maps induce isomorphisms  
\begin{multline*}
   \Symb_{\Gamma_0 (Np), C} (\mD^{\dagger, (k, \epsilon)})^{< h} \isom \Symb_{\Gamma_0 (Np), C} (\mD^{(k, \epsilon)}_r)^{< h} \\ \isom
\Symb_{\Gamma_0 (Np), C} (\mD^{\dagger, (k, \epsilon)}_r)^{< h} \isom \Symb_{\Gamma_0 (Np), C} (\mD^{(k, \epsilon)}_s)^{< h},   
\end{multline*}
which are equivariant for the Hecke operators $U_p$ and $T_v$ for $v \nmid Np$. Consequently, for any $m \ge 1$ and $r> 0$ we have natural isomorphisms 
\[
{\empty} \: \:  \Symb_{\Gamma_0 (Np), C} (\mD^{\dagger, (k, \epsilon)})^{< h} \isom \Symb_{\Gamma_0 (Np^m), C} (\mD^{\dagger, (k, \epsilon)})^{< h} \isom \Symb_{\Gamma_0 (Np^m), C} (\mD_r^{(k, \epsilon)})^{< h}, 
\]
again equivariant for the Hecke operators $U_p$ and $T_v$ for $v \nmid Np$.
\ep 
\bproof 
For the first claim, the same proof of \cite[Lemma 5.3]{Pollack_Stevens_critical} applies in this slightly more general setup, using $U_p$-operators to extend the radius of convergence. The second claim follows from the first and Lemma \ref{lemma:level_raising_and_lowering_at_p}. 
\eproof 

In the next proposition, let $l \ne p$ be a prime satisfying $l \equiv 11 \!\mod 12$ as in Proposition \ref{prop:partial_free_for_l_square} and Corollary \ref{cor:free_over_Gamma}. 

\bp \label{prop:SES_for_Theta_k} 
Let $N = l^2$ and $C = \P^1(\Q) - \Gamma_0 (N) \cdot\{\infty\}$. 
Let $\epsilon$ be a character of conductor $m$.  Write $\Gamma_m := \Gamma_0 (N) \cap \Gamma_1 (p^m)$. Then there is an exact sequence 
\begin{multline*}
\quad    0 \ra \Symb_{\Gamma_0 (Np^m), C} (\mD^{\dagger, (-k-2, \epsilon)}) (\tu{det}^{k+1}) \\ \xrightarrow{ \: \Theta_k \: } \Symb_{\Gamma_0 (Np^m), C} (\mD^{\dagger, (k, \epsilon)}) \lra \Symb_{\Gamma_m, C} (V^k)_{\epsilon} \ra 0, \quad 
\end{multline*}
equivariant for the action by the Hecke operators $U_p$ and $T_v$ for $v \nmid Np$. Here $(-)(\tu{det}^{k+1})$ means that the Hecke action of $[\Gamma_m s \Gamma_m]$ is modified by multiplying the factor $\det(s)^{k+1}$. 
\ep 

\bproof 
By   \cite[Proposition 3.3]{Pollack_Stevens_critical} and its proof, we have a short exact sequence 
\[
0 \lra \mD^{\dagger, (-k-2, \epsilon)} \xrightarrow{\;\Theta_k\;} \mD^{\dagger, (k, \epsilon)} \lra V^{(k, \epsilon)} \lra 0 
\]
of $\Gamma_0 (Np^m)$-modules, where $\Theta_k$ is induced from the map 
$(d/dz)^{k+1}: \mA^{\dagger, (k, \epsilon)} \ra \mA^{\dagger, (-2-k, \epsilon)}$ of $\Gamma_0 (N p^m)$-modules\footnote{Because of the existence of $\epsilon$, the map $(d/dz)^{k+1}$ is not compatible with the $\Sigma_0(p)$-action in general.} 
and $V^{(k, \epsilon)}$ is defined in Remark \ref{remark:epsilon_component}. The assertion then follows from Corollary \ref{cor:free_over_Gamma} and Remark \ref{remark:epsilon_component}. 
\eproof 

\br Although we will not need it, let us remark that by the same proof of \cite[Lemma 5.1]{Pollack_Stevens_critical}, Proposition \ref{prop:SES_for_Theta_k} holds for more general choices of $C$ and more general tame level. This uses the cohomological interpretation of partial modular symbols from \cite[Proposition 2.4]{BD} (see also \cite[Proposition 4.1]{BD}). 
\er

From Proposition \ref{prop:SES_for_Theta_k}, we deduce the following slightly more refined version of Stevens's control theorem (which is the modular symbol analogue of \cite[Proposition 4]{Buzzard_auto_family}) 
\bc[Stevens's control theorem] \label{cor:control_thm}
 Retain assumptions from Proposition \ref{prop:SES_for_Theta_k}. 
The specialization map (\ref{eq:specialize_to_classical}) induces an isomorphism
\[
  \Symb_{\Gamma_0 (Np), C} (\mD^{\dagger, (k, \epsilon)})^{<k+1} \isom \Symb_{\Gamma_0(N) \cap \Gamma_1(p^m), C} (V^k)_{\epsilon}^{<k+1},
\]
which is equivariant for the Hecke operators $U_p$ and $T_v$ for $v \nmid Np$.
\ec 

\bproof 
By the same  argument of \cite[Theorem 5.4]{Pollack_Stevens_critical}, we have 
\[
\big(\Symb_{\Gamma_0 (Np^m), C} (\mD^{\dagger, (-k-2, \epsilon)}) (\tu{det}^{k+1}) \big)^{< (k+1)} = 0. 
\]
The claim then follows from Proposition \ref{prop:finite_slope_witness_no_p_level}, Proposition \ref{prop:SES_for_Theta_k}, and \cite[Proposition 4.8]{Pollack_Stevens_critical}.
\eproof

\subsection{Integral models and characteristic power series} \label{ss:integral_models} In this subsection we consider $\Lambda$ coefficients. Let us denote the ideal $(p, T)  $ inside $ \Lambda$ by $\fm_{\Lambda}$. 
For any topological ring $A$ where $p$ is topologically nilpotent, let $\mC(\Z_p; A)$ denote the $A$-module of continuous functions $\Z_p \ra A$. Recall that $\mC(\Z_p; \Z_p)$ admits an orthornormal basis 
\begin{equation}\label{eq:Mahler_basis}
    \textstyle 1, z, z^{[2]}:= {\binom{z}{2}}, ...,  z^{[n]}:= {\binom{z}{n}}, ... 
\end{equation}
over $\Z_p$ with respect to the supremum norm $|f| := \max_{z \in \Z_p} |f(x)|$, which in turn gives rise to a basis for  $\mC(\Z_p; \Lambda)$ over $\Lambda$. The basis in (\ref{eq:Mahler_basis}) is called the \emph{Mahler basis} of the space $\mC(\Z_p; \Z_p)$ (resp. $\mC(\Z_p; \Lambda)$). 
We refer the reader to \cite{LWX, Colmez} for expositions on Mahler expansions and related $p$-adic analysis. Next let $\chi_{\Lambda}: \Z_p^\times \ra \Lambda^\times$ be the universal character. Equip $\mC (\Z_p, \Lambda)$ with an action of $\Sigma_0 (p)$ given by 
\begin{equation} \label{eq:action_of_Gamma_on_C(Z_p,Lambda)}
    (\gamma \cdot \psi) (z) = \chi_{\Lambda} (a+cz) \psi (\frac{b+dz}{a+cz}). 
\end{equation}
Let $\mC(\Z_p; \Lambda)^{\vee}$ denote the 
linear dual of $\mC(\Z_p; \Lambda)$  and let $\mD^{\tu{int}}$ denote its submodule consisting of convergent sums of the form $\sum_{i \ge 0} c_i \delta^{[i]}, c_i \in \Lambda$, where  $\delta^{[i]} := (z^{[i]})^{\vee}$ denotes the dual Mahler basis and $\Lambda$ is equipped with the $\fm_{\Lambda}$-adic topology.
The action (\ref{eq:action_of_Gamma_on_C(Z_p,Lambda)}) induces an action of $\Sigma_0 (p)$ on $\mC (\Z_p; \Lambda)^\vee$ which preserves $\mD^{\tu{int}}.$ Later we will see that partial modular symbols valued in $\mD^{\tu{int}}$ provides a suitable integral model of overconvergent partial modular symbols. 

\br \label{remark:basis_suffices} For later use we introduce a variant of $D_r(L)$ discussed in Subsection \ref{ss:oc_partial_symbols}. First note that $D_r(L)$ can be identified with the set of distributions $\sum d_n \delta^{[n]}$ with coefficients $d_n \in L$ such that $|d_n|p^{-nr/(p-1)}$ is bounded as $n \ra \infty$ (see \cite[Theorem 1.4.7]{Colmez} and \cite[Section 1.3]{JN}). We define $\sq D_r(L) \subset D_r(L)$ to be the subspace that consists of elements of the form $\sum d_n \delta^{[n]}$ such that $|d_n|p^{-nr/(p-1)} \ra 0$ as $n \ra \infty$.\footnote{Our $D_r$ is denoted by $\mD^{< p^{-r/(p-1)}}$ in \cite{JN}, and our $\sq D_r$ is denoted by $\mD^{p^{-r/(p-1)}}$ there.} The space $\sq D_r(L)$ has orthonormal basis $\{p^{c_n}\delta^{[n]}\}_{n \ge 0}$ where $c_n = - \lfloor \frac{nr}{p-1} \rfloor$ (see \cite[Section 3.1, in particular, 3.2.1]{JN}). Taking the inverse limit, we have 
\[\mD^{\dagger}_{r_0}(L) =  \varprojlim_{r > r_0} 
D_r(L) \cong  \varprojlim_{r > r_0} 
\sq D_r(L)
\]
for every $r_0$ with $0 \le r_0 < 1$. 
\er

\subsection{The Atkin--Lehner involution} \label{ss:AL}
Keep Assumption \ref{assumption:sec3} as before.
It would be important for our applications that a certain Atkin--Lehner action preserves the cusps $C$. To this end, let us choose a matrix  $\big(\begin{smallmatrix} a & b \\ Nc & p^md \end{smallmatrix} \big) \in \SL_2(\Z)$ where $c, d \in \Z$. 

\bl 
Let $C = \P^1(\Q) - \Gamma_0 (N) \cdot \{\infty\}$. 
Then the action of the matrix 
\[
W = \begin{pmatrix} a p^m & b \\ Np^m c & p^m d \end{pmatrix}
\]
preserves $C$. 
\el 

\bproof 
Suppose that $W x = \gamma \infty$ for some $\gamma =  \big(\begin{smallmatrix} a' & b' \\ N c'& d' \end{smallmatrix} \big) \in \Gamma_0 (N)$, then $x$ is of the form
\[
x = \frac{p^m d a' - N b c'}{p^m N (ac'-a'c)} \in \Gamma_0 (N) \cdot \{\infty\}. 
\]
\eproof

Consequently, for a right $\Sigma_W$-module $V$ where $\Sigma_W$ denotes the submonoid of $\GL_2(\Q)$ generated by $\Sigma_0 (p)$ and $W$, we obtain a well-defined Atkin--Lehner involution on the partial modular symbols $\Symb_{\Gamma, C} (V)$. Let us specialize to our case of interest. Let $N = l^2$ as in Proposition \ref{prop:partial_free_for_l_square} and further suppose that $\Gamma_0 (N)$ contains no $2$- or $3$- torsion elements. Let $e_1, ..., e_t$ be a choice of generators of $\Delta_C^0$ as a free $\Gamma_0(N)$-module. Let $L$ be a finite extension of $\Q_p$ and let $(k, \epsilon)$ be a classical weight character of conductor $m$, then 
\[
\Symb_{\Gamma_0(N) \cap \Gamma_1(p^m), C} (V^k)_{\epsilon} \cong \Symb_{\Gamma_0 (Np^m), C} (V^{(k, \epsilon)})
\]
is an $L$-vector space of dimension $d_{k,m}:= p^{m-1}(p+1)t(k+1)$. The following assertion is a variant of \cite[Proposition 3.22]{LWX}. 

\bp \label{prop:Atkin_Lehner}
Let $a_0 \le a_1 \le ... \le  a_{d_{k, m}-1}$ (resp. $b_0 \le b_1 \le ... \le b_{d_{k, m}-1 } $) be the slopes of the $U_p$ action on $
\Symb_{\Gamma_0(N p^m), C} (V^{k,\epsilon})$ (resp. $\Symb_{\Gamma_0(N p^m), C} (V^{k,\epsilon^{-1}})$), then we have $a_i + b_i = k+1$ for each $i \in \{0, 1, ..., d_{k, m}-1\}$.  
\ep 

\bproof 
We claim that, after possibly replacing $L$ by a finite extension, the inclusion 
\[
\Symb_{\Gamma_0(N)\cap \Gamma_1(p^m), C} (V^k)^{\tu{s.s.}} \subset M_{k+2}(\Gamma_0 (N) \cap \Gamma_1(p^m), L)^{\oplus 2}
\]
from Part (2) of Proposition \ref{prop:Eichler_Shimura} is compatible with the Atkin--Lehner involution induced by $W$. The proposition then follows from this using the same argument as \cite[Proposition 3.22]{LWX}. Now let us justify the claim. For simplicity we write $\Gamma = \Gamma_0 (N) \cap \Gamma_1(p^m)$. From the proof of  Proposition \ref{prop:Eichler_Shimura}, it suffices to show that the isomorphisms in (\ref{eq:hecke_equivariant_ES}) are compatible with the Atkin--Lehner action. Upon fixing an isomorphism between $\C_p$ and $\C$, it suffices to consider the isomorphisms in 
 (\ref{eq:hecke_equivariant_ES}) over $\C$. To this end, let us first observe that the  map 
\[f \mapsto \phi_f:  
S_{k+2} (\Gamma, \C) \lra \Symb_{\Gamma} (V^k(\R))
\]
defined by the formula 
\[ 
\phi_f(\{x\} - \{y\}) (P(z)) := \tu{Re} \int_{y}^{x} f(z) P(z)dz
\]
is compatible with the double coset action by $[\Gamma s \Gamma]$ for all $s \in \GL_2(\Q)$ with $\det s > 0$, and in particular with the Atkin--Lehner action (which is identified with $[\Gamma W \Gamma]$). After complexifying, this in turn implies that  (by the proof of \cite[Proposition 2.5]{BD}) the first isomorphism in (\ref{eq:hecke_equivariant_ES}) over $\C$ is compatible with the Atkin--Lehner involution, as desired. The compatibility for the other isomorphism can be deduced by a similar argument. 
\eproof

\subsection{The eigencurve of partial modular symbols} \label{ss:eigencurve_of_modular_symbols} 
In this subsection we recall the eigencurve of partial modular symbols from \cite[Section 4.4]{BD}. 

Keep Assumption \ref{assumption:sec3} as before. Take $\Gamma = \Gamma_0(Np)$ and suppose that the assumptions of Remark \ref{remark:Hecke_action} and Remark \ref{remark:diamond} are satisfied.
Let $\iota$ denote the involution on partial modular symbols induced from the double coset action by $[\Gamma \big(\begin{smallmatrix} 1 & 0 \\ 0 & -1 \end{smallmatrix} \big) \Gamma]$. For a $\Q_p$-Banach algebra $L$, this involution splits $\Symb_{\Gamma, C}(D_r(L))$ into two eigenspaces $\Symb_{\Gamma, C}(D_r(L))^{\pm}$ for $\iota$. Let us fix a sign $\pm$. For each affinoid subdomain $X = \tu{Sp} A \subset \mW$, the $U_p$ operator acts compactly on $\Symb_{\Gamma, C} (D_r(A))$, where the $\Sigma_0(p)$-action on $D_r(A)$  
is defined by (\ref{eq:action_on_D_r}) and (\ref{eq:action_on_A_r}) using the canonical (or tautological) character 
\begin{equation} \label{eq:tautological_char}
    \Z_p^\times \xrightarrow{\chi_{\Lambda}} \Lambda^\times \ra A^\times.
\end{equation} We then apply Buzzard's eigenvariety machine (\cite[Construction 5.7]{Buzzard_eigen}) to obtain the eigencurve $\mX_{\tu{symb}, C}^{\pm}$ of partial modular symbols supported on $C$ with tame level $N$ and sign $\pm$, where we take $U_p$ to be Buzzard's operator $\phi$ (the requirements to apply \cite{Buzzard_eigen} can be checked the same way as \cite{BD}. See Proposition 4.3 of \emph{loc.cit.}, also see \cite{eigenbook}). The eigencurve $\mX_{\tu{symb}, C}^{\pm}$ enjoys similar properties as the usual eigencurve $\mX$ constructed by Coleman--Mazur and Buzzard for modular forms. In particular, it is a reduced rigid analytic curve over $\Q_p$ equipped with a locally finite flat weight map $\tu{wt}: \mX_{\tu{symb}, C}^{\pm} \ra \mW$ and a slope map $a_p: \mX_{\tu{symb}, C}^{\pm} \ra \G_m^{\tu{rig}}$. Moreover, we have the following comparison result. 

\bl \label{lemma:eigencurve_contains_cuspidal_part}
Fix the tame level $\Gamma_0 (N)$ and sign $\pm$ as above.  There are (unique) closed immersions 
\[
\mX^{\tu{cusp}} \hookrightarrow \mX_{\tu{symb}, C}^{\pm} \hookrightarrow \mX 
\]
compatible with the weight map and the slope map, where $\mX^{\tu{cusp}}$ denotes the cuspidal part of $\mX$. In particular, $
\mX^{\tu{cusp}}, \mX_{\tu{symb}, C}^{\pm}$, and $\mX$ agree away from the ordinary component. 
\el 

\bproof 
This follows from Chenevier's comparison theorem (\cite[Proposition 3.2]{Chenevier}) and Proposition \ref{prop:Eichler_Shimura}. 
\eproof

\section{Estimates of the Newton Polygon} 

In this section, we analyze the Newton polygon of the $U_p$ operator in terms of the Mahler basis.  Recall that we have fixed an isomorphism $\Z_p^\times \cong \F_p^\times \times (1+p\Z_p)$ from the beginning of the introduction, which induces an isomorphism 
\[
\Lambda \cong \Z_p[\F_p^\times] \otimes \Z_p [\![T]\!]
\]
where $T$ corresponds to the element $[\exp(p)]-1 \in \Lambda$. As in \cite{LWX}, we let 
\[\Lambda^{> 1/p}:= \Lambda [\![pT^{-1}]\!] \cong \Z_p[\F_p^\times] \otimes \Z_p[\![T, pT^{-1}]\!].\]
The following lemma is essentially a reformulation of \cite[Proposition 3.14]{LWX}. 
\bl \label{lemma:estimate_on_delta_action}
\be 
\item Let $\gamma =  \big(\begin{smallmatrix} a & b \\ c & d \end{smallmatrix} \big) \in \big(\begin{smallmatrix} \Z_p^\times & \Z_p \\ p \Z_p & p \Z_p  \end{smallmatrix} \big) $ and let $P (\gamma) = (P_{m,n}(\gamma))$ denote the (infinite) matrix for the right action of $\gamma$ on $\mD^{\tu{int}}$ with respect to the basis $\{\delta^{[i]}\}$ from Section \ref{ss:integral_models}. Then the entries of $P(\gamma)$ satisfies 
\[
P_{m, n} (\gamma) \in \fm_{\Lambda}^{\max (n - \lfloor m/p \rfloor, 0)},
\]
where $\fm_{\Lambda}$ denotes the ideal $(p, T)  $ inside $ \Lambda$. 
\item Similarly, when $\gamma \in \Sigma_0 (p)$, we have 
\[
P_{m, n} (\gamma) \in \fm_{\Lambda}^{\max (n-m, 0)}.
\]
\ee 
\el 
\bproof 
It suffices to show that the matrix $P'$ for the left action\footnote{
Note that the convention for the action of $\gamma$ on $\mC(\Z_p; \Lambda)$ we use is different from \cite{LWX}. In our setup, the action is a left action given by (\ref{eq:action_of_Gamma_on_C(Z_p,Lambda)})
(compare with the right action from \cite[Equation 2.3.2]{LWX}). This explains why the shape of matrix in Part (1) (resp. our definition of $\Sigma_0 (p)$) is different from the one used in \textit{loc.cit.} (resp. $\mathbf{M}_1$ in \textit{loc.cit.}). Also see Remark \ref{remark:transpose}.}
of $\gamma$ on $\mC(\Z_p; \Lambda)$ with respect to the Mahler basis $\{z^{[i]}\}$ satisfies $P'_{m, n} \in \fm_{\Lambda}^{\max (m - \lfloor n/p \rfloor, 0)}$ for Part (1) (resp. $P'_{m, n} \in \fm_{\Lambda}^{\max (m - n, 0)}$ for Part (2)). For this the proof of \cite[Proposition 3.14]{LWX} carries over \emph{verbatim}. 
\eproof

For the next proposition and its corollary let us assume that $\Gamma =\Gamma_0 (pN)$ where $N = l^2$ satisfies the assumptions in Proposition \ref{prop:partial_free_for_l_square}, and let $C = \P^1(\Q) - \Gamma_0 (N) \cdot \{\infty\}$. Then $\Gamma_0 (pN) \subset \Gamma_0 (N)$ has index $s = p+1$, where a coset representative can be chosen to be $\{\eta_{0}, ..., \eta_{p-1}, \eta_{\infty}\}$ where $\eta_i = \big(\begin{smallmatrix} 1 & 0 \\ Ni & 1 \end{smallmatrix} \big)$ for $i = 0, ..., p-1$ and $\eta_\infty = \big(\begin{smallmatrix} a & b \\ N & p \end{smallmatrix} \big) $ for some choice of $a, b \in \Z$ such that $\det (\eta_\infty) = 1$. By Proposition \ref{prop:partial_free_for_l_square}, $\Delta^0_C$ is a free $\Z[\Gamma]$-module of rank $st = (p+1)t$. Let $\sq e_1, ..., \sq e_{st}$ be a set of generators, which then 
determines an isomorphism of $\Lambda$-modules
\begin{equation} \label{eq:free_over_C_Lambda}
\rho: \Symb_{\Gamma, C} (\mD^{\tu{int}}) \isom  (\mD^{\tu{int}})^{\oplus st} 
\end{equation}
sending $\varphi$ to $(\varphi(\sq e_i))$ (also see Corollary \ref{cor:free_over_Gamma}).

\bp  \label{prop:char_power_series_Up}
Let $P = (P_{m, n})_{m, n \in \Z_{\ge 0}}$ denote the infinite matrix of the $U_p$ operator using the isomorphism  (\ref{eq:free_over_C_Lambda}), with respect to the following basis 
\begin{equation} \label{eq:Mahler_basis_t_copies}
\delta^{[0]}_0, \dots, \delta^{[0]}_{st-1}, \delta^{[1]}_0, \dots, \delta^{[1]}_{st-1}, \delta^{[2]}_0, \dots, \delta^{[2]}_{st-1} , \dots    
\end{equation}
  where $\delta_i^{[0]}, \delta_i^{[1]}, ... \delta_i^{[n]}, ... $ denote the basis of the $i^{th}$ copy of $\mD^{\tu{int}}$ from (\ref{eq:free_over_C_Lambda}). 
Then 
\be
\item 
The characteristic power series
\begin{equation} \label{eq:char_power_series}
F(X) 
:= \varinjlim_{n \ra \infty} \det (1 - (P_{i, j})_{i, j \le n} X) \in \Lambda [\![X]\!]
\end{equation} 
is well-defined. 
\item Fix an affinoid subdomain $\tu{Sp} A \subset \mW$. For each $r \in (0, 1)$, 
write $\mD_r^{\dagger, \chi_{\Lambda}} = \mD_r^{\dagger, \chi_{\Lambda}}(A)$ for $\mD_r^{\dagger} (A)$ with a $\Sigma_0 (p)$-action given by (\ref{eq:action_on_A_r}) and  (\ref{eq:action_on_D_r}) using the character (\ref{eq:tautological_char}). Then the power series $F(X)$ defined in (\ref{eq:char_power_series}) from Part (1)\footnote{More precisely, the image of $F(X)$ along the map $\Lambda [\![X]\!] \ra A[\![X]\!]$ induced by the structure map $\Lambda \ra A$.} agrees with the characteristic power series of the $U_p$ operator acting on $\Symb_{\Gamma, C} (\mD_r^{\dagger, \chi_{\Lambda}}).$ 
\item Write $F(X) = \sum_{n \ge 0} c_n X^n \in \Lambda [\![X]\!]$ for the characteristic power series in (1), then   
\[
c_n \in T^{\lambda(n)} \cdot \Lambda^{>1/p} 
\]
for all $n \ge 0$, where $\lambda (0):= 0$ and $\lambda(n)$ is defined by the recurrence relation \[
\lambda(n) - \lambda (n-1) = \lfloor \frac{n-1}{(p+1)t} \rfloor - \lfloor \frac{n-1}{p(p+1)t} \rfloor. 
\]
\ee 
\ep 

 
\bproof 
The proof is similar to that of \cite[Theorem 3.16]{LWX}. Let us first establish an analogue of \cite[Proposition 3.1 (1) and (3)]{LWX}. 
For each $\varphi \in \Symb_{\Gamma, C} (\mD^{\tu{int}})$ we have $\rho (\varphi) = (\varphi(\sq e_i))$ via the isomorphism (\ref{eq:free_over_C_Lambda}). Write $\gamma_a =  \big(\begin{smallmatrix} 1 & a \\ 0 & p \end{smallmatrix} \big)$, now we have
\begin{align*}
    U_p \cdot \varphi (\sq e_i) &  = \:  \sum_{a = 0}^{p-1} \varphi (\gamma_a \cdot \sq e_i) |\!|_{\gamma_a} \\
    & = \: \sum_{a} \varphi(\sum_{l \in I_a} \gamma_l^{a, i} \cdot \sq e_{a_l}) |\!|_{\gamma_a} \\
    & = \: \sum_{a} \sum_{l \in I_a}  \varphi(\sq e_{a_l}) |\!|_{(\gamma_l^{a, i})^{-1} \cdot \gamma_a}, 
\end{align*}
where, for the second equality, we use the fact the modular path $\gamma_a \cdot \sq e_i$ can be written as a finite sum of unimodular paths (see Subsection \ref{ss:Delta_free_over_Gamma}), each of which is of the form  $\gamma_l^{a, i} \cdot \sq e_{a_l}$ for some $\gamma_l^{a, i} \in \Gamma$ and $a_l \in \{0, ..., st-1\}$ where $s = p+1$ (the $a_l$'s are allowed to be repetitive). Note that since each $\gamma_l^{a, i} \in \Gamma \subset \Gamma_0 (p)$, we have 
\[
(\gamma_l^{a, i})^{-1} \cdot \gamma_a \in \begin{pmatrix} 
\Z_p^\times & \Z_p \\ p \Z_p & p \Z_p  
\end{pmatrix}
\]
when we view $(\gamma_l^{a, i})^{-1} \cdot \gamma_a$ as an element in $M_2(\Z_p)$. Therefore, via the isomorphism (\ref{eq:free_over_C_Lambda}), the $U_p$ operator can be represented by a $st \times st$ matrix for each element in $\bigoplus_{0}^{st-1} \mD^{\tu{int}}$, where each entry of the $st \times st$ matrix is a finite sum of operators of the form $(-)|\!|_{\delta}$ with $\delta \in  \big(\begin{smallmatrix} \Z_p^\times & \Z_p \\ p \Z_p & p \Z_p  \end{smallmatrix} \big).$ The rest of the argument is almost identical to that of \cite[Theorem 3.16]{LWX}, which we briefly explain for completeness. First note that, in terms of the basis in (\ref{eq:Mahler_basis_t_copies}), we know that the matrix $P$ for the $U_p$ operator satisfies 
\begin{equation}\label{eq:Lambda_estimate_P_mn}
P_{n,m} \in \fm_{\Lambda}^{\max (\lfloor m/st \rfloor - \lfloor n/pst \rfloor, 0)}
\end{equation}
by Lemma \ref{lemma:estimate_on_delta_action}. This implies that for any sufficiently large integer $j$, the matrix $P$ is strictly lower  triangular modulo $\fm_{\Lambda}^j$ except for the first $\lfloor \frac{pstj}{p-1} \rfloor \times \lfloor \frac{pstj}{p-1} \rfloor$-minor, so $F(X) \in \Lambda [\![X]\!]$ is indeed well-defined. 

For the  claim in Part (2), by Remark \ref{remark:basis_suffices} it suffices to show that $F(X)$ agrees with the characteristic power series of the $U_p$ operator acting on each $\Symb_{\Gamma, C} (\sq D_{r}^{\chi_r})$. Again from Remark \ref{remark:basis_suffices}, we know that
\begin{equation} \label{eq:scaling_basis_matrix}
p^{c_0}\delta^{[0]}_0, \dots, p^{c_0}\delta^{[0]}_{st-1}, p^{c_1}\delta^{[1]}_0, \dots, p^{c_1}\delta^{[1]}_{st-1}, p^{c_2}\delta^{[2]}_0, \dots, p^{c_2}\delta^{[2]}_{st-1} , \dots    
\end{equation} 
forms an orthonormal basis for $\Symb_{\Gamma, C} (\sq D_{r}^{\chi_r})$, thus the matrix for the $U_p$ operator with respect to (\ref{eq:scaling_basis_matrix}) differs from $P$ by conjugating by an (infinite) diagonal matrix. The claim then follows from this (see also Proposition 2.17 of \emph{loc.cit.}).

Part (3) of the proposition follows from the same proof of Theorem 3.16 of \emph{loc.cit}. More precisely, we conjugate the matrix $P$ by the diagonal matrix $\mathcal T$ with entries \[(1,..., 1, T, ..., T, T^2, ..., T^2, ... )\] with $t$-copies of each $T^i$ to form $P' = \mathcal T  P \mathcal T^{-1}$, then for all $m, n \ge 0$ we have
\[
P_{n, m}' \in T^{\lfloor n/st \rfloor - \lfloor n/pst \rfloor} \Lambda^{> 1/p} 
\]
using the estimate in  (\ref{eq:Lambda_estimate_P_mn}) and that $p = T \cdot p/T \in \Lambda^{> 1/p}$. In other words, the entries of the $n^{th}$-row of $P'$ all lie in  $ T^{\lfloor n/st \rfloor - \lfloor n/pst \rfloor} \Lambda^{> 1/p}$.  This proves the proposition.
\eproof 

\bd \label{def:Newt_polygon_over_weight}
Now fix a character $\omega$ of $\F_p^\times$, which induces a quotient map $\tu{pr}_\omega: \Lambda \ra \Z_p [\![T]\!]$ by evaluating $\F_p^\times$ using $\omega$. We then define 
\[
F_\omega (X) = \sum_{n \ge 0} c_{\omega, n} X^n
\]
to be the image of the characteristic power series $F(X)$ of $U_p$ from (\ref{eq:char_power_series}) along the projection $\tu{pr}_\omega$. Note that each $c_{\omega, n}$ is an element in $\Z_p[\![T]\!]$ of the form 
\begin{equation} \label{eq:def_of_b_n}
c_{\omega, n} = c_{\omega, n}(T) = \sum_{m \ge 0} b_{\omega, n, m} T^m.
\end{equation} For each $\beta \in \C_p$ with $|\beta| < 1$, we denote the \emph{specialization} of $F_\omega(X)$ to $\beta$ by \[F_{\omega}^{\beta} (X) = \sum_n c_{\omega, n}(\beta) X^n \in \C_p[\![X]\!],\] and denote the Newton polygon of $F_\omega^{\beta}(X)$ by $\tu{Newt}_{\omega, \beta}$. 
\ed 

Following the convention of \cite{LWX}, we call the polygon 
\begin{equation} \label{eq:LB_polygon}
\tu{LB}_{\beta} = \{(n, \lambda(n)\cdot v(\beta))\}
\end{equation}
  the \emph{lower bound polygon} at $\beta$. 
Proposition \ref{prop:char_power_series_Up} implies the following important corollary. 
\bc \label{cor:Newt_vs_LB} 
Let $\beta \in \C_p$ be an element with valuation $v(\beta) \in (0, 1)$. Let $b_{\omega, n, m}$ be as defined in (\ref{eq:def_of_b_n}). Then for each $\omega$ and $n$, we have 
\[
v(b_{\omega,n, m} \beta^m) \ge \max \{\lambda(n) - m, 0\} + m \cdot v(\beta) 
\]
for every $m \ge 0$. Consequently, we have
\[
\begin{cases} 
v(c_{\omega, n} (\beta)) = \lambda (n) \cdot v(\beta) & \tu{ if (and only if) }  b_{\omega, n, \lambda (n)} \in \Z_p^\times \\
v(c_{\omega, n} (\beta)) \ge \lambda (n) \cdot v(\beta) + \min(v(\beta), 1 - v(\beta)) & \tu{ if }  b_{\omega, n, \lambda (n)} \notin  \Z_p^\times, \end{cases}
\] where $b_{\omega, n, \lambda (n)}$ are defined in (\ref{eq:def_of_b_n}). In particular, the Newton polygon $\tu{Newt}_{\omega, \beta}$ of $F_{\omega}^{\beta}(X)$ 
lies on or above the lower bound polygon $\tu{LB}_{\beta}$ at each $\beta$ with $v(\beta) \in (0, 1)$. 
\ec

\bproof 
This follows from the same argument of \cite[Corollary 3.18]{LWX}, using Proposition \ref{prop:char_power_series_Up} Part (2). 
\eproof

\section{Proof of the main theorem}

Let us fix a component $\omega$ of $\mW$. 
Consider a classical weight $\chi = (k,\epsilon)$ of conductor $2$ with $\chi|_{\F_p^\times} = \omega$. Recall that the  $T$-coordinate $T_{(k, \epsilon)} \in \C_p$ of $(k, \epsilon)$ satisfies $v(T_{(k, \epsilon)}) \in (0, 1)$.  Following \cite{LWX}, we let $n_k = p k st = p(p+1) k t$. 

\bl \label{lemma:touching_of_polygon}
For each $k \ge 1$ and $p$-Nebentypus character $\epsilon$ of conductor $2$, the sum of the first $n_{k}$ $U_p$-slopes in $\Symb_{\Gamma_0 (N p), C} (\mD^{\dagger, (k-1, \epsilon)})$ is $pk^2st = p(p+1)k^2 t$. Moreover, the Newton polygon $\tu{Newt}_{\omega, T_{(k-1, \epsilon)}}$ of $F_{\omega}^{T_{(k-1, \epsilon)}} (X) =  \sum_{n} c_{\omega, n}(T_{(k-1, \epsilon)}) X^n$ goes through the point
\[(n_{k}, \lambda(n_{k}) v(T_{(k-1, \epsilon)})) = (p(p+1)kt,  p(p+1)k^2t/2) \in \tu{LB}_{T_{(k-1, \epsilon)}}.\]
\el 
 
\bproof 
This is Step 1 of the proof of \cite[Theorem 1.3]{LWX} (see Section 3.23 of \emph{loc.cit}). We need to use Proposition \ref{prop:Atkin_Lehner} and Corollary \ref{cor:control_thm} in place of Proposition 3.22 and 2.15 of \emph{loc.cit.} 
\eproof 

Lemma \ref{lemma:touching_of_polygon} essentially says that at each classical point $(k-1, \epsilon)$ in the weight space, we know precisely one point on the Newton polygon (specialized to $T_{(k-1, \epsilon)}$ in the sense of Definition \ref{def:Newt_polygon_over_weight}). Corollary \ref{cor:Newt_vs_LB} allows us to ``propagate the touching of the polygons across the boundary  $\mW^{>\frac{1}{p}}$ of the weight space", which gives the first part of the Halo conjecture. More precisely, 

\bt \label{mainthm:disjoint_union}
Fix a character $\omega$ and let $\mX_{\omega}$ denote the eigencurve over the $\omega$-component of the weight space $\mW$. Then the space $\mX_{\omega}^{>\frac{1}{p}}$ breaks into a disjoint union 
\begin{equation} \label{eq:initial_decomposition}
\mX_{\omega}^{>\frac{1}{p}} = \mX_{\omega, 0} \sqcup \mX_{\omega, (0, 1)} \sqcup \mX_{\omega, 1} \sqcup \mX_{\omega, (1, 2)} \sqcup...
\end{equation}
of rigid analytic spaces, each of which is finite flat over $\mW_{\omega}^{>\frac{1}{p}}$ via the weight map, such that, for each point $x \in \mX_{\omega, I}$ where $I$ denotes the interval $[n, n]$ or $(n, n+1)$, its slope satisfies 
\[
v (a_p(x)) \in (p-1) v(T_{\tu{wt}(x)}) \cdot I
\]  
\et 

\bproof 
It suffices to prove an analogous statement for the spectral curve and the same argument of \cite{LWX} applies in this setup. For completeness, let us recall the key part of the argument. It is easy to compute the endpoints of the line segment in $\tu{LB}_\beta$ containing the point $(n_k, \lambda(n_k) v(\beta))$: they have $x$-coordinates $n_k -st$ and $n_k + st$ respectively. Now let $n_k^{-}$ (resp. $n_k^{+}$) denote the $x$-coordinate of the left (resp. right) endpoint of the line segment in  $\tu{Newt}_{\omega, T_{(k-1, \epsilon)}}$ which contains the point $(n_k, pk^2st/2) = (n_k, p(p+1)k^2t/2)$.\footnote{Note that we  have $n_k - st \le n_k^{-} \le n_k \le n_k^{+} \le n_k + st$. By convention, we set $n_k^{-} = n_k^+ = n_k$ if $(n_k, pk^2st/2)$ is a vertex of $\tu{Newt}_{\omega, T_{(k-1, \epsilon)}}$.} Moreover, this entire line segment  (between $n_k^{-}$ and $n_k^{+}$) has slope $k(p-1)v(T_{(k-1,\epsilon)}) = k$ and belongs to $\tu{LB}_{T_{(k-1, \epsilon)}}$. 
The latter fact has the following crucial consequence: for every weight $x \in \mW_{\omega}^{>1/p}$, the points 
\[
(n_k^{-}, \lambda(n_k^{-}) \cdot v(T_x)), \tu{ and } (n_k^{+}, \lambda(n_k^{+}) \cdot v(T_x))
\]
are the two vertices of a line segment in $\tu{Newt}_{\omega, T_{x}}$.\footnote{They could be the same vertex, in the case when $n_k^{-}= n_k^{+}$.} This is because, by Corollary \ref{cor:Newt_vs_LB}, we can interpret $n_k^{-}$ as the minimal index $j \in [n_{k}-st, n_k]$ such that $b_{\omega, j, \lambda(j)} \in \Z_p^\times$ is a $p$-adic unit, which is a condition that is independent of $x \in \mW_{\omega}^{> 1/p}$. The same argument applies to $n_k^{+}$. The decomposition of the spectral curve follows from this, and finite flatness follows from \cite[Corollary 4.3]{Buzzard_eigen}. 
\eproof 

The proof of Theorem \ref{mainthm:disjoint_union} above in fact implies the following. Recall that $n_k = pkst$. 
\bc 
For each $\beta \in \C_p$ with $v(\beta) \in (0, 1)$, we define the \textit{upper bound polygon} at $\beta$ to be
\begin{equation}
    \tu{UB}_\beta = \{ (n_k, \lambda(n_k) v(\beta))\}. 
\end{equation}
Then the Newton polygon  $\tu{Newt}_{\omega, \beta}$ at $\beta$ lies on or below $\tu{UB}_{\beta}$ for each $\beta$ with $v(\beta) \in (0, 1)$. 
\ec

\br[Remark on the degree] 
Let $\mY_{\omega}$ (resp. $\mY_{\omega}^{>1/p}$) denote the spectral curve over $\mW_{\omega}$ (resp. the pullback of $\mY_{\omega}$ to $\mW_{\omega}^{> 1/p}$), by the proof above we know that $\mY_{\omega}^{>1/p}$ admits a decomposition $\mY_{\omega}^{>1/p} = \bigsqcup_{I} \mY_{\omega, I}$ similar to (\ref{eq:initial_decomposition}), where $I$ runs through intervals of form $[n,n]$ and $(n, n+1)$ for all $n \in \N$.  The degree of each $\mY_{\omega, I}$ can be determined similarly as in \cite{LWX}. More precisely, we have 
\begin{align*}
\deg \mY_{\omega, k} \qquad \; &= \begin{cases} 
d_{\tu{ord}} (\omega) & \tu{ if } k = 0 \\ 
\frac{1}{2} \Big(r_{\tu{ord}, C} (\omega^{-1} \omega_0^{2k -2}) + r_{\tu{ord}, C} (\omega \omega_0^{-2k}) \Big) & \tu{ if } k \ge 1
\end{cases} \\
\deg \mY_{\omega, (k, k+1)} & = 
\frac{1}{2} \Big(pt - r_{\tu{ord}, C} (\omega^{-1} \omega_0^{2k}) - r_{\tu{ord}, C} (\omega \omega_0^{-2k}) \Big). 
\end{align*}
Here $d_{\tu{ord}}(\omega)$ (resp. $r_{\tu{ord}, C}(\omega)$) denotes the dimension of the ordinary subspace of of modular forms (resp. of partial modular symbols supported on $C$) of weight $2$ and character $\omega$, and $\omega_0$ denotes the character $\F_p^\times \cong (\Z_p^\times)_{\tu{tor}} \hookrightarrow \Z_p^\times$ where the second map is the natural inclusion. Note that this makes use of Proposition \ref{prop:SES_for_Theta_k} and Proposition \ref{prop:Eichler_Shimura} (which explains the factor $\frac{1}{2}$ in the formula, see also Subsection \ref{ss:eigencurve_of_modular_symbols}, in particular, Lemma \ref{lemma:eigencurve_contains_cuspidal_part}). 
\er 
 
Going further into the boundary, we are  ready to establish the full Halo conjecture (Conjecture \ref{conj:Halo}), which asserts that sufficiently close to the boundary (i.e., when $v(\beta)$ is sufficiently small), the ``shape'' of the Newton polygon $\tu{Newt}_{\omega, \beta}$ is independent of $\beta$. To formulate the precise statement, let us set 
\[\lambda = p^{-8/(st(p^2-1)+8)} = p^{-8/((p-1)(p+1)^2t+8)}.\]
\br 
We follow \cite{LWX} for the choice of $\lambda$, which is chosen so that for every $\beta \in \C_p$ such that $0 <v(\beta) < \lambda$, the point $(n, \lambda (n) v(\beta)+1)$ lies strictly above the polygon $\tu{UB}_{\beta}$ (see Lemma 4.1 of \emph{loc.cit}.). 
\er  

The following Proposition is extracted from the proof of \cite[Theorem 1.5]{LWX}. We briefly explain the argument for completeness. 
\bp[\cite{LWX}] \label{prop:key_estimate}
If $(n, v(c_{\omega, n}(\beta_0)))$ is a point lying strictly below the upper bound polygon $\tu{UB}_{\beta_0}$ for some $\beta_0$ with $v(\beta_0) \in (0, \lambda)$, then there exists an integer $h_n \ge \lambda(n)$ such that for every $\beta \in \C_p$ with $v(\beta) \in (0, \lambda)$, we have $v(c_{\omega, n}(\beta)) = h_n \cdot v(\beta)$. 
\ep

\bproof 
Let $h_n = h_n (\beta_0)$ be the minimal index $m$ such that $ v (c_{\omega, n}(\beta_0)) \ge v (b_{\omega, n, m} \beta_0^m)$ (which \emph{a priori} depends on $\beta_0$). We first claim that, if $(n, y)$ is a point strictly below $\tu{UB}_{\beta}$ with $v(\beta) \in (0, \lambda)$, then 
\begin{equation} \label{eq:m_smaller_than_lambda_case}
y < \lambda(n) v(\beta) + \frac{(p-1)^2st \cdot v(\beta)}{8} < v(b_{\omega, n, m} \beta^m) \quad \textup{ for all } m < \lambda(n).
\end{equation}
Here the two inequalities follow from \cite[Lemma 4.1]{LWX} and Corollary \ref{cor:Newt_vs_LB} respectively (see also \cite[Inequality 4.2.1 \& 4.2.2]{LWX}).  In particular, we know that $h_n \ge \lambda (n)$ by the assumptions. Moreover, the choice of $\lambda$ forces $b_{\omega, n, h_n} \in \Z_p^\times$ to be a unit, which in turn implies that 
\begin{equation} \label{eq:m_bigger_than_h_case}
v(b_{\omega, n, m} \beta^m) > v (b_{\omega, n, h_n} \beta^{h_n}) \quad  
\end{equation} for all $m > h_n$ and all $\beta$ with $v(\beta) \in (0, \lambda)$. We want to show that, 
\begin{equation}\label{eq:m_is_h_n}
v(c_{\omega, n} (\beta)) = v (b_{\omega, n, h_n} \beta^{h_n}) = h_n \cdot v(\beta)
\end{equation}
for every $\beta$ with $v(\beta) \in (0, \lambda)$. Now observe that for $\beta = \beta_0$, the equality (\ref{eq:m_is_h_n}) follows from the minimality of $h_n$. This then implies  that, for each $\beta$ we have $(n, h_n \cdot v(\beta))$ lying strictly below $\tu{UB}_{\beta}$. Thus, by (\ref{eq:m_smaller_than_lambda_case}), we have 
\begin{equation} \label{eq:smaller_case}
    v (b_{\omega, n, h_n} \beta^{h_n}) <  v (b_{\omega, n, m} \beta^{m}) \quad 
\end{equation}   for all $m < \lambda(n)$. In fact, the equality (\ref{eq:m_is_h_n}) for $\beta = \beta_0$ also implies that $b_{\omega, n, m} \in p\Z_p$ when $m < h_n$. Therefore, 
\begin{equation} \label{eq:the_case_m_between_h_and_lambda}
   v (b_{\omega, n, m} \beta^{m}) \ge 1 + \lambda(n) v(\beta) > h_n \cdot v(\beta) \quad   
\end{equation}
if $\lambda(n) \le m < h_n$. Here the second equality uses the choice of $\lambda.$ Together, (\ref{eq:smaller_case}), (\ref{eq:the_case_m_between_h_and_lambda}) and (\ref{eq:m_bigger_than_h_case}) give the desired (\ref{eq:m_is_h_n}). 
\eproof 

\bt \label{mainthm:arithmic_progression} 
Fix a character $\omega: \F_p^\times \ra \Z_p^\times$. Then there exists a sequence of increasing rational numbers  
$\alpha_0, \alpha_1, \alpha_2, \dots $ 
(depending on $\omega$) tending to $\infty$, such that $\mX^{> \lambda}_{\omega}$ breaks into a disjoint union \[
\mX^{> \lambda}_{\omega} = \mX_{\omega, \alpha_0 } \sqcup \mX_{\omega, \alpha_1 } \sqcup \mX_{\omega, \alpha_2 } \sqcup \dots 
\]
with the property that for any point $x \in \mX_{\omega, \alpha_i }$, we have 
\[
v(a_p(x)) = (p-1) \alpha_i   \cdot  v(T_{\tu{wt}(x)}). 
\]
Moreover, these $\alpha_i  $'s (counted with multiplicity $\deg \mX_{\omega, \alpha_i}$) form a finite union of arithmetic progressions.
\et 

\bproof 
Proposition \ref{prop:key_estimate} implies that there exists a subset $I \subset \N$ such that, for all $\beta$ with $v(\beta) \in (0, \lambda)$, the Newton polygon $\tu{Newt}_{\omega, \beta}$ is the convex hull of the points
\[
\left\{\big(n_k, \lambda(n_k) \cdot v(\beta) \big) \right\}_{k \in \N}, \quad \left\{\big(n, h_n \cdot v(\beta) \big) \right\}_{n \in I}.
\]
Note that $\tu{Newt}_{\beta}$ only depends on $v(\beta)$. Moreover, the $x$-coordinates of the breakpoints of these Newton polygons are independent of $\beta$. The existence of $\alpha_i$'s and the decomposition of $\mX_{\omega}^{> \lambda}$ are immediate from this. This also implies the assertion on arithmetic progressions, using the argument of the proof of \cite[Theorem 1.5]{LWX}, or equivalently, the proof of \cite[Theorem B]{Bergdall_Pollack}. 
\eproof

\bibliographystyle{plain}
\bibliography{halo}

\vspace{0.7cm} 
\end{document}